\documentclass{amsart}
\usepackage[utf8]{inputenc}
\usepackage{amsmath,amstext, amsthm, empheq, extpfeil,  
amsbsy,amssymb,marvosym,fancyhdr,graphicx,amscd,amsfonts,latexsym,delarray,stackrel,lineno,color,cite,appendix,xcolor,url,hyperref,slashed,enumerate,fullpage}
\usepackage{wasysym}
\usepackage{subfig}

\usepackage{srcltx}
\usepackage{mathrsfs}
\usepackage{float} 

\usepackage[all]{xy}
\usepackage{t1enc}
\usepackage{mathrsfs}
\usepackage{pifont}

\usepackage{mathpple}
\usepackage[T1]{fontenc}

\definecolor{Green}{rgb}{0,1,0}
\definecolor{Blue}{RGB}{0,0,191}
\definecolor{mathmodecolor}{RGB}{0,102,0}
\definecolor{keywordcolor}{RGB}{0,51,151}
\definecolor{sourcebackgroundcolor}{RGB}{255,247,223}
\definecolor{unixagred}{RGB}{255,0,0}
\definecolor{lightgray}{RGB}{191,191,191}
\definecolor{green}{RGB}{1,191,191}

\newcommand*\patchAmsMathEnvironmentForLineno[1]{%
  \expandafter\let\csname old#1\expandafter\endcsname\csname #1\endcsname
  \expandafter\let\csname oldend#1\expandafter\endcsname\csname end#1\endcsname
  \renewenvironment{#1}%
     {\linenomath\csname old#1\endcsname}%
     {\csname oldend#1\endcsname\endlinenomath}}%
\newcommand*\patchBothAmsMathEnvironmentsForLineno[1]{%
  \patchAmsMathEnvironmentForLineno{#1}%
  \patchAmsMathEnvironmentForLineno{#1*}}%
\AtBeginDocument{%
\patchBothAmsMathEnvironmentsForLineno{equation}%
\patchBothAmsMathEnvironmentsForLineno{align}%
\patchBothAmsMathEnvironmentsForLineno{flalign}%
\patchBothAmsMathEnvironmentsForLineno{alignat}%
\patchBothAmsMathEnvironmentsForLineno{gather}%
\patchBothAmsMathEnvironmentsForLineno{multline}%
}

\newtheorem{thm}{Theorem}[section]
\newtheorem{prop}[thm]{Proposition}
\newtheorem{cor}[thm]{Corollary}
\newtheorem{lem}[thm]{Lemma}

\newtheorem{rem}[thm]{Remark}

\newtheorem{fact}[thm]{Fact}

\def\Det{{\rm Det}}

\def\id{{\rm id}}

\def\tr{{\rm Tr}}

\def\A{{\mathbb A}}
\def\B{{\mathbb B}}
\def\C{{\mathbb C}}
\def\F{{\mathbb F}}

\def\N{{\mathbb N}}

\def\Q{{\mathbb Q}}
\def\R{{\mathbb R}}
\def\Z{{\mathbb Z}}

\def\tr{{\rm tr}}

\def\cD{{\mathcal D}}

\def\cH{{\mathcal H}}

\def\cS{{\mathcal S}}

\newcommand{\ie}{{\it i.e.\/}\ }

\def\sin{{{\rm sin}}}
\def\cos{{{\rm cos}}}

\def\id{{\mbox{Id}}}

 \def\scal2{{\mathscr S}}
 
\def\fourier{\F}
\def\sr0{{\cS^{\rm ev}_0}}

\def\sar0{{\cS_0(\A_\Q)}}

\def\fourier{\F}

\def\F{{\mathbb F}}

\def\bm2{{\rm \B mod^2}}
\def\b2{{\rm \B mod^{\mathfrak s}}}

\def\id{{\rm id}}

\def\tr{{\rm tr}}

\def\A{{\mathbb A}}
\def\C{{\mathbb C}}
\def\F{{\mathbb F}}

\def\N{{\mathbb N}}

\def\Q{{\mathbb Q}}
\def\R{{\mathbb R}}

\def\Z{{\mathbb Z}}

\def\B{{\mathbb B}}

\def\fourier{\F}

\def\tr{{\rm tr}}

\def\cD{{\mathcal D}}

\def\cH{{\mathcal H}}

\def\cS{{\mathcal S}}

\def\part{\partial}

\def\sin{{{\rm sin}}}
\def\cos{{{\rm cos}}}

\def\id{{\mbox{Id}}}

\definecolor{trust}{rgb}{0,1,1}

\def\Det{{\rm Det}}

\def\tilde{\widetilde}

\newcommand{\nil}[1]{}

 \title{Quadratic Forms, Real Zeros \\[2ex] and Echoes of the Spectral Action}
\author{Alain Connes}
\address{College de France, 3 rue Ulm, F75005, Paris, France\\
I.H.E.S. F-91440 Bures-sur-Yvette, France}
\email{alain@connes.org}

\author{Walter D. van Suijlekom}

\address{Institute for Mathematics, Astrophysics and Particle Physics, Radboud
University Nijmegen, Heyendaalseweg 135, 6525 AJ Nijmegen, The Netherlands.
}
\email{waltervs@math.ru.nl}

\begin{document}
\maketitle
\begin{center}
  \large\itshape 
  Dedicated to Huzihiro Araki\\[0.5em]
  with gratitude and admiration
\end{center}
\begin{abstract}

	For a  real distribution  $\mathcal{D}$  on the interval $[0,L]$ with $\tilde{\mathcal{ D}}$ the associated even distribution on the interval $[-L, L]$, we prove that if the associated quadratic form with Schwartz kernel $\tilde{\mathcal{D}}(x - y)$ defines a lower-bounded selfadjoint operator  on $L^2([-\frac{L}{2}, \frac{L}{2}])$, whose lowest spectral value $\lambda$ is a simple, isolated eigenvalue with even eigenfunction $\xi$, then all the zeros of the entire function $\widehat \xi(z)$, the Fourier transform of $\xi$, lie on the real line.

The proof proceeds in five steps.
(1) We give a $C^*$-algebraic proof of a corollary of Carathéodory--Fej\'er’s 1911 structure theorem for Toeplitz matrices: if $T \in M_n(\mathbb{C})$ is a Hermitian, positive semidefinite Toeplitz matrix of rank $n - 1$, and $\xi \in \ker T$, then the polynomial $P(z) = \sum \xi_j z^j$ has all its zeros on the unit circle.
(2) We formulate and prove a continuous analogue of this result, replacing the Toeplitz matrix with a convolution operator with continuous kernel $h(x - y)$, and the polynomial $P(z)$ with the Fourier transform of the eigenfunction corresponding to the largest eigenvalue.
(3) We analyze finite-dimensional truncations of the quadratic forms defined by real, even distributions $\mathcal{D}$ on $[-L, L]$, and observe that the resulting matrices exhibit a structure previously encountered in perturbative expansions of the spectral action.
(4) We establish an analogue of Carathéodory--Fej\'er’s corollary for matrices of this specific structure, thereby extending the zero localization result beyond the classical Toeplitz setting.
(5) Finally, we apply a classical theorem of Hurwitz concerning the zeros of uniform limits of holomorphic functions to deduce the general result stated above.
\end{abstract}
\tableofcontents
\section{Introduction}
The Carathéodory-Fej\'er theorem from 1911~\cite{CF11} (see also \cite[Theorem 1.3.6]{BW}) describes the structure of Hermitian, positive semidefinite Toeplitz matrices as follows. Let

\[
T = \begin{bmatrix}
c_0 & \overline{c}_1 & \overline{c}_2 & \cdots & \overline{c}_{n} \\
c_1 & c_0 & \overline{c}_1 & \cdots & \overline{c}_{n-1} \\
c_2 & c_1 & c_0 & \cdots & \overline{c}_{n-2} \\
\vdots & \vdots & \vdots & \ddots & \vdots \\
c_{n} & c_{n-1} & c_{n-2} & \cdots & c_0
\end{bmatrix}.
\]

Then, if \(T\) is positive semidefinite of rank \(r\), there exist distinct points \( z_1, \dots, z_r \in \mathbb{T} \subset \mathbb{C} \) (the unit circle), and positive weights \( \alpha_1, \dots, \alpha_r > 0 \), such that
\[
T = V D V^*,
\]
where \( V \in \mathbb{C}^{(n+1) \times r} \) is the Vandermonde matrix
\[
V = \begin{bmatrix}
1 & 1 & \cdots & 1 \\
z_1 & z_2 & \cdots & z_r \\
z_1^2 & z_2^2 & \cdots & z_r^2 \\
\vdots & \vdots & \ddots & \vdots \\
z_1^{n} & z_2^{n} & \cdots & z_r^{n}
\end{bmatrix},
\]
and \( D = \mathrm{diag}(\alpha_1, \dots, \alpha_r) \) is a diagonal matrix with positive real entries.

In our work on operator systems~\cite{CW}, we showed how the above theorem can be derived from the duality theory of operator systems. This factorization also plays a central role in the proof of truncated Weil positivity in~\cite{ckw}.

A direct corollary of the Carathéodory-Fej\'er theorem is the following:

\begin{cor}\label{corcar}
Let \( T \in M_{n+1}(\mathbb{C}) \) be a Hermitian, positive semidefinite Toeplitz matrix of rank \( n  \), and let \( \xi \in \ker T \). Then all the zeros of the polynomial
\[
P(z) := \sum_{j=0}^n \xi_{j} z^j
\]
lie on the unit circle.
\end{cor}

This corollary exhibits a striking number-theoretic flavor, resonating with the analogue of the Riemann Hypothesis for function fields; see~\cite{HP} for a further discussion of this connection. In number theory, Toeplitz matrices of this kind naturally arise, and Corollary~\ref{corcar} applies to show that the zeros of the polynomial \( P(z) \), associated to an eigenvector for the smallest eigenvalue of such a matrix, all lie on the unit circle. The key difficulty in this context, then, becomes the verification that zero is indeed the (simple) minimal eigenvalue of \(T\).

In the present paper, we investigate a distributional analogue of Corollary~\ref{corcar}, motivated by its potential relevance to the Riemann Hypothesis itself. Our approach proceeds through several stages:

\begin{enumerate}
  \item We give a $C^*$-algebraic proof of Corollary~\ref{corcar}.
  \item We formulate and prove a continuous analogue of Corollary~\ref{corcar}, in which the Toeplitz matrix is replaced by a convolution operator with continuous kernel \( h(x - y) \), and the polynomial \( P(z) \) by the Fourier transform of the eigenfunction corresponding to the largest eigenvalue.
  \item We analyze the finite truncations of quadratic forms defined by real even distributions \( \mathcal{D} \) supported on \( [-L, L] \), and observe that the resulting matrices exhibit a structure previously encountered in perturbative expansions of the spectral action.
  \item We prove an analogue of Corollary~\ref{corcar} for matrices of this special type.
  \item Finally, using a classical result of Hurwitz on the zeros of uniform limits of holomorphic functions, we deduce the following general theorem:
  
\begin{thm}\label{mainintro}
Let \( L > 0 \), \( \mathcal{D} \) be a real  distribution on the interval \( [0, L] \) and $\tilde \cD$ the associated even distribution on \( [-L, L] \). Assume that the quadratic form with Schwartz kernel \( \tilde \cD(x - y) \) defines a lower-bounded selfadjoint operator \( A \) on \( L^2([-\frac{L}{2}, \frac{L}{2}]) \), and that the minimum of its spectrum is a simple, isolated eigenvalue \( \lambda \), with even eigenfunction $\xi$. Then all the zeros of the entire function $\widehat \xi(z)$, $z\in \C$, Fourier transform  of $\xi$ lie on the real line.
\end{thm}
\end{enumerate}
We refer to Theorem \ref{main} for the precise formulation, the above formulation is slightly unprecise as shown in Remark \ref{nuance}.
In the course of the proof, we encounter a number of illustrative examples and special cases. A detailed matrix-based verification of the theorem is given in an appendix.

\subsection*{Acknowledgements}
WvS thanks Teun van Nuland for useful discussions. The authors are indebted to an anonymous referee for numerous suggestions and comments. 

\subsection*{Conflict of interest statement}
The authors declare that there is no conflict of interest.

\subsection*{Data availability statement}
Not applicable: no data has been created or analysed in this study.

\section{Toeplitz case}
Recall the operator system $C^*(\Z)_{(n+1)} \subseteq C^*(\Z)$ given by Fourier truncations on the interval $[-n,n] \subset \Z$ from \cite{CW}
It is the dual operator system of the operator system of Toeplitz matrices, which in particular allows to associate a positive linear form $\mathcal L_T$ to any positive Toeplitz matrix $T$:
\begin{align*}
  \mathcal L_T: C^*(\Z)_{(n+1)} &\to \C; \qquad 
f \mapsto \sum_{k=-n}^n f_{k} c_k.
  \end{align*}
Note that a  positive real Toeplitz matrix $T$ of size $n+1$ can be written in the following form:
$$
T = \begin{pmatrix}
  c_0 & c_1 & \ldots & c_n\\
  c_1 & c_0 & \ddots &\vdots \\
  \vdots &\ddots  & \ddots &  c_1 \\
  c_n & \ldots & c_1 & c_0 
    \end{pmatrix} ; \qquad (c_k \in \R).
$$
The starting point for our proof of Corollary \ref{corcar} is the following purely $*$-algebraic result:
\begin{prop}
  Let $T = (c_k)$ be a positive $(n+1)$-dimensional real Toeplitz matrix of rank $n$ and let $(a_j)_{j=0}^n$  be a real vector in $\ker T$.
  \begin{enumerate}
  \item The ideal $J$ of $A=\C[X,X^{-1}]$ generated by $P=\sum_0^n a_j X^j$ is stable under the canonical involution, $(aX^n)^*:=\overline{a}X^{-n}$, of $\C[X,X^{-1}]$.
    \item The monomials $X^j$ ($j = 0,\ldots, n-1$) form a basis of $A/J$.
  \item There exists a unique linear form $\phi$ on the quotient $A/J$ such that 
  \begin{equation}\label{phi}
  \phi(X^j)=c_j, \forall j\in \{0,\ldots,n-1\}
 \end{equation}
    \item The linear form $\phi$ is positive on the $*$-algebra $A/J$.	
    \end{enumerate}
\end{prop}

\proof
(1)
We have that $J^* = J$ since $P(X)$ is either palendromic, or anti-palendromic, {\em i.e.} $a_{n-j} = \pm a_j$ for all $j =0, \ldots, n$. Indeed, because of the structure of $T$ as a Toeplitz matrix, it follows that if $(a_j)_j$ is in the kernel of $T$, so is $(a_{n-j})_j$. Since furthermore this kernel is one-dimensional, it follows that $a_{n-j} = \lambda a_j$, hence $a_j = \lambda^2 a_j$, which implies $\lambda = \pm 1$. But then
$$
P^*(X) = \sum_{j=0}^ n {a_j} X^{-j} =  \pm X ^{-n} P(X)
$$
as claimed. 

(2) Since $T$ is of rank $n$ one has $a_0\neq 0$ (see \cite[Lemma 33]{HP} or show directly that if $a_0=0$ then $(a_{j+1})_{j=0}^n$ with $a_{n+1}\equiv 0$ is also in $\ker T$). Thus $X$ is invertible in $\C[X]/J'$ where $J'=P\C[X]$ since modulo $P$ one has $a_0+ X Q=0$ for $Q=(P-a_0)/X$. The ring $A=\C[X,X^{-1}]$ is the localisation of $\C[X]$ at $X$ and localisation commutes with quotients \cite[Proposition 3.3 and Corollary 3.4]{AM}, since $X$ is not a zero divisor modulo $J'$. Since $X$ is invertible in  $\C[X]/J'$ localization at $X$ does not change  $\C[X]/J'$ and one thus gets the equality $A/J=\C[X]/J'$ and the claim follows.


(3) Follows from (2).\newline
(4) Let $ \phi : A \to \C$ be the unique linear form which vanishes on the ideal $J$ and fulfills \eqref{phi}.
In order to show that  $ \phi$ is positive, we first show that $\phi(X^{-j} ) = c_j$ for all $j = 1 ,\ldots, n-1$.
We have
$$
 \phi(X^{-1} P) = 0 \implies a_0  \phi (X^{-1} ) + \sum_{j=1}^n a_j  \phi(X^{j-1}) = 0 \implies a_0  \phi (X^{-1} ) + \sum_{j=1}^n a_j c_{j-1} = 0.
$$
Since $(a_j) \in \ker T$ we have in terms of the second row of $T$ that
$$
\sum_{j=0}^n c_{j-1} a_j = 0 
$$
and hence $\phi(X^{-1}) = c_{-1} = c_1$. This argument can be repeated by considering subsequent rows in $T$ to obtain by induction that $\phi(X^{-k}) = c_{-k} = c_k$. Note that from the first row of $T$ it also follows that $\phi(X^n) = c_n$. We conclude that $\phi(X^j) = c_{|j|}$ for all $j = - n, \ldots, n$.

To show that $\phi$ is positive, \ie that $\phi(f^**f)\geq 0$ for all $f\in A$, note that the value of $\phi(f^**g)$ only depends on the classes of $f$ and $g$ in $A/J$. Thus one can take $f=\sum_0^{n-1} f_jX^j$, one then gets 
$$
\phi(f^**f)=\sum f_j\overline{ f_k}\phi(X^{j-k})=\sum f_j\overline{ f_k}c_{\vert j-k\vert}=\langle f\mid Tf\rangle \geq 0
$$ 
which shows (4).
\endproof

\begin{prop}
\label{toeplitzcase}
  Let $T = (c_k)$ be a positive $(n+1)$-dimensional Toeplitz matrix of rank $n$. If $(a_j)_{j=0}^n$ is a vector in $\ker T$ then the polynomial $P(z) = \sum_j a_j z^j$ has all zeros on unit circle in $\C$. 
\end{prop}
\proof
The positive linear form $\phi$ on $A$ vanishing on $J$ defines a positive linear form on the envelopping $C^*$-algebra $C^*(\Z)$  of the involutive algebra $A$,  \ie a positive measure on the Pontrjagin dual $U(1)$ of $\Z$. 
This measure is supported by the $n$ eigenvalues of the unitary $\pi(X)$ associated to $X$ in the GNS representation $\pi$ of $(A,\phi)$ which is of dimension $n$ by construction.
 Since $P\in J\subset \ker \pi $, these eigenvalues correspond  to the roots of the generator $P$ of $J$ which are hence all of modulus one. 
\endproof

\begin{rem}
In general, when the kernel of $T$ is more than one-dimensional ---in other words, when the extreme eigenvalue is not simple--- it follows from the Carath\'eodory-Fej\'er decomposition that there is an equivalence between being in the kernel and the vanishing of the corresponding polynomial on the complex numbers of modulus one that appear in this decomposition. So if it happens that the number of these complex numbers is strictly less than n, then this condition will be fulfilled by polynomials which will have these particular complex numbers as zeros, but which otherwise can have arbitrary other zeros. This means that it is not true in general if the eigenvalue is not simple, that the theorem holds. The correct formulation of the theorem is that if you take the intersection of the zeros of the various eigenfunctions, then they are all on the unit circle. This is reminiscent to the notion of the radical of a quadratic form. 
  \end{rem}

\section{The continuous kernel case}
In this section we shall extend the Toeplitz case to the continuous case.

\begin{figure}
  \centering
  \includegraphics[scale=0.4]{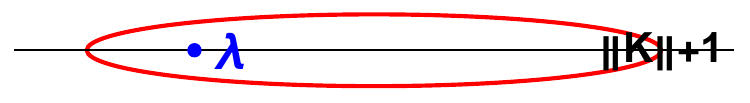}
  \caption{The contour of integration for the definition of the projection $P$ in Equation \eqref{eq:proj}.}
  \label{fig:contour}
\end{figure}

\begin{thm}\label{thmtoepcont} 
Let $h\in C([-L,L])$ be an even real continuous function. Let $K$ be the operator on $L^2([0,L])$
given by 
\begin{equation}\label{opk}
(Kf)(x)=\int h(x-y)f(y)dy	.
\end{equation}
Then $K$ is a compact selfadjoint operator. Assume that its eigenvalue of largest modulus is simple and let $\xi\in L^2([0,L])$ be its eigenvector. If we extend $\xi$ to an element of $L^2(\R)$ to be zero outside $[0,L]$ then all the zeros of the entire function $\widehat \xi$ belong to $\R\subset \C$.
\end{thm}
\proof The operator $K$ of \eqref{opk} is of Hilbert-Schmidt class since $h(x-y)$ is square integrable. We first approximate $K$ by finite rank operators as follows. Since the function $h$ is real continuous and even on $[-L,L]$ we can find, given $\epsilon>0$ an $\alpha >0$ and scalars $h_j \in \R$ for $ j=0,\ldots, L/\alpha=N\in \N$, such that, with $j(x)$ denoting the integer part of $x/\alpha$, one has 
$$
\vert h(x-y)-h_{\vert j(x)-j(y)\vert}\vert \leq \epsilon, \qquad \forall x,y \in [0,L].
$$
This follows provided one chooses $\alpha >0$ and the scalars $h_j \in \R$ such that $\forall x\in [0,L]$
$$
\vert h(x)-h_{j}\vert \leq \epsilon, \qquad \forall j \mid \vert j-j(x)\vert \leq 1, 
$$
as one gets by comparing $j(\vert x-y\vert)$ with $\vert j(x)-j(y)\vert$.
\begin{figure}
  \centering
\includegraphics[scale=0.3]{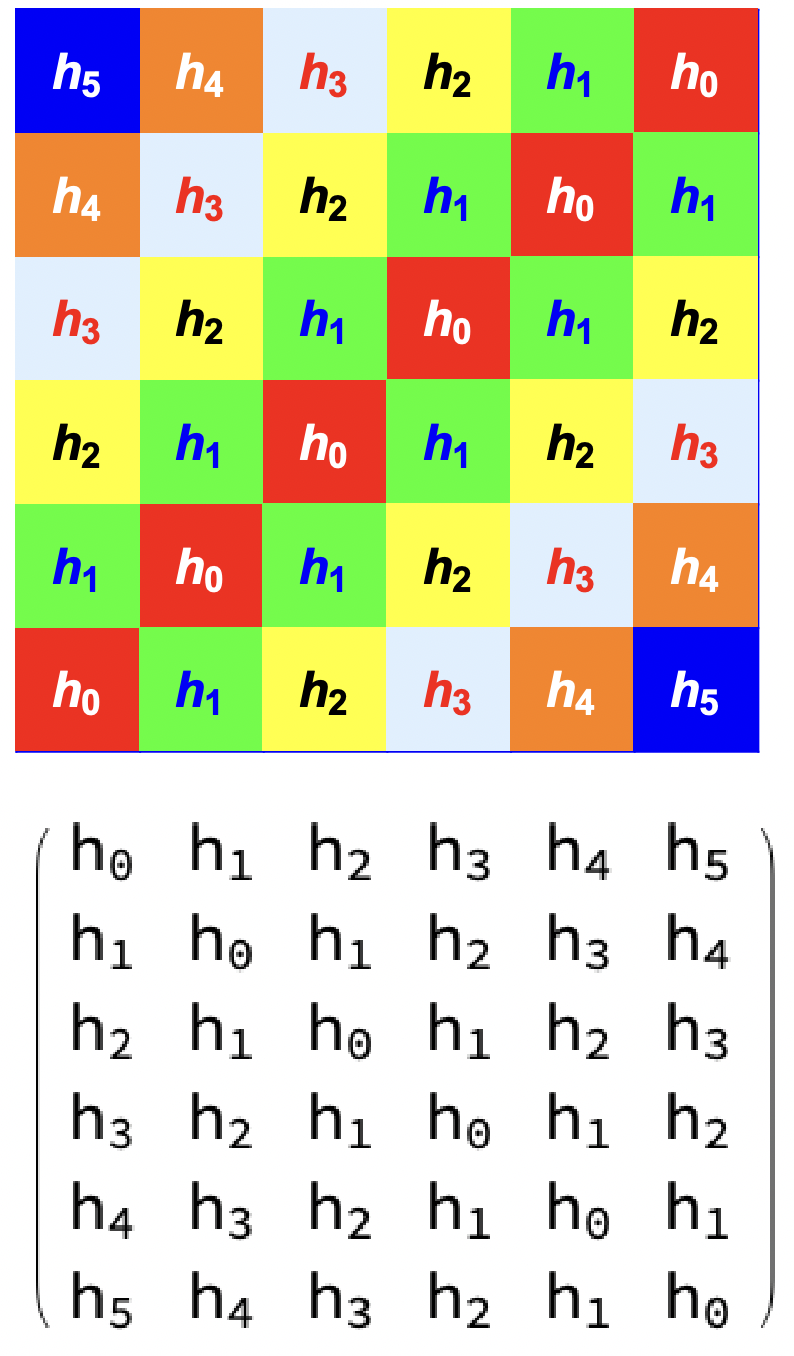}\label{toep}
\caption{The approximation of $h(\vert x-y\vert)$ and the reflected Toeplitz matrix (by symmetry with respect to the $x$ axis).}
\end{figure}
Let then $\chi_j$ be the characteristic function of the interval $I_j:=\{x\mid j(x)=j\}$ and $T$ the Toeplitz matrix 
$
T_{n,m}:= h_{\vert n-m\vert}
$, one thus obtains the inequality 
$$
\vert h(x-y)-\sum T_{n,m} \chi_n(x)\chi_m(y)\vert \leq \epsilon, \qquad  \forall x,y \in [0,L].
$$
It follows from the Hilbert-Schmidt control of the norm that one obtains in this manner a sequence $R_n$ of finite rank operators converging in norm to $K$, each of the form, with $T$ a real symmetric Toeplitz matrix and $\chi_j$ characteristic functions of intervals,
$$
R_n=\sum T_{i,j} |\chi_i \rangle \langle \chi_j|, \  \  \Vert K-R_n\Vert \to 0.
$$ 
Let then $\xi$ be an eigenvector of norm $1$ for the eigenvalue $\lambda$ of $K$ of largest modulus; without loss of generality we assume that $\lambda > 0$. The spectral projection $P$ obtained as a Cauchy integral of the resolvent of $K$  along the contour $C$ (see Figure \ref{fig:contour}), isolating $\lambda$ from the rest of ${\rm Spec}K$,
\begin{equation}
  P=\frac{1}{2\pi i}\int_C(z-K)^{-1}dz
  \label{eq:proj}
\end{equation}
is of rank one and fulfills $P\xi=\xi$. Let $P_n$ be defined by the same formula using $R_n$ in place of $K$. This makes sense for $n$ large enough since $\Vert K-R_n\Vert\to 0$ when $n\to \infty$. One has $P_n\to P$ in norm and thus $P_n$ is of rank one for $n$ large enough. Let $\xi_n=P_n\xi$. One has $\xi_n\to \xi$ in norm. Let us show that for $n$ large enough, $\xi_n$ is an eigenvector of $R_n$ for the largest eigenvalue.
By construction, for $n$ large enough, it is a non-zero eigenvector of $R_n$ and is the only one for the interval between the lower part of $C$ and its highest part. Since we have chosen this highest part to be $\Vert K\Vert +1$, it follows that, for $n$ large enough, it is the highest eigenvalue of $R_n$. We can thus assert that $\xi$ is the norm limit of the sequence $\xi_n$ where, for each $n$, $\xi_n$ is of the form,   
$$
\xi_n=\sum a_j \chi_j, \ \ \sum a_j z^j=0\ \Rightarrow \vert z\vert =1
$$
and where the $\chi_j$ are the characteristic functions of $N=L/\alpha$ intervals $I_j=I_0+j\,\alpha$, $0\leq j<N$ forming a partition of $[0,L)$. We now compute the Fourier transform of $\xi_n$. The Fourier transform of $\chi_0$ is 
  $$\widehat \chi_0(s)=\int_0^{\alpha } \exp (-i s x) \, dx=\frac{i \left(-1+e^{-i \alpha  s}\right)}{s}.
  $$
One has $\chi_j(x)=\chi_0(x-j\,\alpha)$, and this gives $\widehat \chi_j(s)=\exp(-is\alpha j)\widehat \chi_0(s)$ so that we obtain
$$
\widehat \xi_n(s)=P(\exp(-is\alpha))\widehat \chi_0(s), \ \ P(z):=\sum a_j z^j
$$
Thus the zeros of $\widehat \xi_n$ are the union of the set $\{\frac{2n\pi}{\alpha}\mid n\in \Z, n\neq 0\}$ of zeros of $\widehat \chi_0$ with the set of complex numbers $s$ such that $\exp(-is\alpha)$ is one of the roots of $P(z)$. But we know (from the Toeplitz case, Proposition \ref{toeplitzcase}) that all these roots are of modulus $1$. For each root $z_k$ of $P(z)=0$,
let $s_k\in \R$ be such that $\exp(-is_k\alpha)=z_k$,  the set $Z$ of zeros of $\widehat \xi_n$ is then
$$
Z=\{\frac{2n\pi}{\alpha}\mid n\in \Z, n\neq 0\}\cup_k \{s_k+\frac{2n\pi}{\alpha}\mid n\in \Z\}\subset \R.
$$
Now since the sequence $\xi_n$ converges to $\xi$ in $L^2[0,L]$, we find for
the entire functions $\widehat \xi_n(s)$ and $\widehat \xi(s)$ that
$$
| \widehat \xi(s) - \widehat \xi_n(s)| \leq \| \xi - \xi_n \|_{L^2} \int_0 ^L e^{2 \Im (s) x}.
$$
We conclude from this that the sequence $\widehat \xi_n(s)$ converges uniformly to $\widehat \xi(s)$ on compact subsets of $\C$. Thus the Hurwitz's theorem shows that all the zeros of  $\widehat \xi(s)$ are real. \endproof

\section{Quadratic form $Q$ associated to a distribution}
Let $L>0$. We start with a  distribution $\cD$ on the interval $[0,L]$, \ie a continuous  linear form 
on $C^\infty([0,L])$ written formally as 
$$
\cD(f)=\int_0^L f(x)\cD(x), \ \ \forall f\in C^\infty([0,L])
$$
and we use it to define an even distribution  $\tilde \cD$ on $[-L,L]$ by symmetrization, 
\begin{equation}\label{symet}
\tilde \cD(f):=\cD(f)+\cD(\tilde f), \ \ \tilde f(x):=f(-x),\  \forall x \in [-L,L].
\end{equation} 
Note that $\tilde \cD(f)$ continues to make sense when $f$ is smooth when restricted to both $[0,L]$ and $[-L,0]$
but not necessarily smooth at $0$. 
 We consider the Hilbert space $\cH:=L^2([0,L],dx)$ with orthonormal basis 
\begin{equation}\label{basis}
U_n(x):=L^{-\frac 12}\exp(2\pi i nx/L), \ \forall x	\in [0,L], \ n\in \Z.
\end{equation}
These functions are extended to $x\in \R$ so that they vanish for $x\notin [0,L]$. One then uses the involution $f^*(x):=\overline{ f(-x)}$ and the convolution 
$$
(f*g)(y):=\int f(x)g(y-x)dx.
$$
We then use $\cD$ to get a densely defined hermitian form $Q$ on trigonometric polynomials by 
\begin{equation}\label{qua0}
\langle f\vert g\rangle_Q:=\tilde\cD(f^**g)=\int_0^L \left((f^**g)(x)+(f^**g)(-x)\right)\cD(x),
\end{equation}
whose matrix in the orthonormal basis $U_n$ is given by
\begin{equation}\label{quad}
\langle U_m\vert U_n\rangle_Q:=\int_0^L\, \left((U_m^**U_n)(y)+(U_m^**U_n)(-y)\right)\cD(y)dy.
\end{equation}
We assume that this defines a lower-bounded self-adjoint quadratic form with spectrum having an isolated eigenvalue at its minimum  and are interested in the eigenvector $\eta$ for this lowest eigenvalue (assumed simple). \newline
For $\cD=\delta_0$ the Dirac mass at $y=0$, one has 
$$
\int_0^L\, \left((U_m^**U_n)(y)+(U_m^**U_n)(-y)\right)\cD(y)dy=2(U_m^**U_n)(0)=2\int_0^LU_n(x)\overline{U_m(x)}dx=2\langle U_m\vert U_n\rangle
$$
which shows that by adding to $\cD$ a multiple of $\delta_0$ one can assume that the quadratic form $Q$ is positive and the eigenvector $\eta$ is in the radical of the quadratic form, \ie that one has 
$$
\langle U_n\vert \eta\rangle_Q=0, \ \ \forall n\in \Z.
$$
By analogy with the Toeplitz case, one would like to obtain a positive linear form on the quotient of the convolution algebra of functions on $\R$  by the ideal generated by $\eta$ allowing one to extend the quadratic form $Q$. This would then yield  an analogue of Corollary \ref{corcar} by showing that  the zeros of the entire function $\widehat \eta$ are all real.

Instead of working directly with this infinite dimensional situation, our strategy is to   first compute the matrix \eqref{quad}, observe that it is of a particular form already met in perturbation theory of the spectral action, and prove a matrix form of the reality of zeros of Fourier transforms of lowest eigenvectors. The infinite dimensional result will then follow by approximation using Hurwitz theorem as in the proof of Theorem \ref{thmtoepcont}.
\subsection{Matrix of the quadratic form $Q$}
For $y\in [0,L]$ one has, for $n\neq m$, using $U_m(t)=0$ for $t<0$ and $U_n(x)=0$ for $x>L$,
$$
(U_m^**U_n)(y)=\int U_m^*(y-x)U_n(x)dx=\int \overline{ U_m(x-y)}U_n(x)dx=\frac 1L\int_y^L\exp(2\pi im(y-x)/L+2\pi inx/L)dx=
$$
$$
\frac{\exp(2\pi im y/L)}{L}\int_y^L\exp(2\pi i(n-m)x/L)dx=\frac{\exp(2\pi im y/L)}{2\pi i(n-m)}\left(\exp(2\pi i(n-m)x/L)\right)_y^L=$$ $$=\frac{\exp(2\pi im y/L)-\exp(2\pi in y/L)}{2\pi i(n-m)}.
$$
Moreover and still with $y\in [0,L]$ one has
$$
(U_m^**U_n)(-y)=\overline{(U_m^**U_n)^*(y)}=\overline{(U_n^**U_m)(y)}.
$$
Thus since the formula $\frac{\exp(2\pi im y/L)-\exp(2\pi in y/L)}{2\pi i(n-m)}$ is symmetric in $n,m$, one obtains for $n\neq m$ and $y\in [0,L]$,
\begin{equation}\label{symreal}
(U_m^**U_n)(y)+(U_m^**U_n)(-y)=2\,\Re\left(\frac{\exp(2\pi im y/L)-\exp(2\pi in y/L)}{2\pi i(n-m)}\right)=	\frac{\sin(2\pi m y/L)-\sin(2\pi n y/L)}{\pi (n-m)}.
\end{equation}
For $m=n$ the same computation gives 
$$
(U_n^**U_n)(y)=\frac 1L\int_y^L\exp(2\pi in(y-x)/L+2\pi inx/L)dx=(1-y/L)\exp(2\pi in y/L)
$$
and
\begin{equation}\label{symreal1}
(U_n^**U_n)(y)+(U_n^**U_n)(-y)=2\,\Re\left((1-y/L)\exp(2\pi in y/L)\right)=2(1-y/L)\cos(2\pi n y/L).\end{equation}
We can summarize the above computation as follows

\begin{prop}\label{matrixcomp} Let $\cD$ be as above and $Q$ be the quadratic form $Q$ of \eqref{quad}. Let $\psi(x):=\frac 1 \pi \int_0^L \sin(2\pi x(1- y/L))\cD(y)dy$. The matrix elements $q_{m,n}$ of $Q$ are given as follows 
\begin{equation}\label{matelem}
	q_{m,n}=\begin{cases} \frac{\psi(m)-\psi(n)}{m-n} & \text{if } n \neq m, \\
		\psi'(n), & \text{if } n = m.
	\end{cases}
\end{equation} 	
\end{prop}
\proof One has, for $n\neq m$,  using \eqref{quad}, \eqref{symreal} and the equality $\sin(2\pi x(1- y/L)=-\sin(2\pi x y/L)$ for $x\in \Z$,
$$
\langle U_m\vert U_n\rangle_Q=\int_0^L\, \left((U_m^**U_n)(y)+(U_m^**U_n)(-y)\right)\cD(y)dy=
$$
$$
\int_0^L \frac{\sin(2\pi m y/L)-\sin(2\pi n y/L)}{\pi (n-m)}\cD(y)dy= \frac{\psi(n)-\psi(m)}{n-m}.
$$
For $n=m$, one has  using \eqref{quad}, \eqref{symreal1}
$$
\langle U_n\vert U_n\rangle_Q=2\int_0^L (1-y/L)\cos(2\pi n y/L)\cD(y)dy=\partial_x\psi(x)\vert_{x=n}
$$
which gives the required equality.
\endproof
\subsection{The diagonal values}
Proposition \ref{matrixcomp} gives the diagonal values of the matrix of the quadratic form $Q$, but unlike the off-diagonal values which only depend upon the first Fourier  components of the distribution $\cD$, the diagonal values involve all the Fourier components of $\cD$. This follows from the equalities 
\begin{align*}
2 \int_0^1 (1-x) \exp (2 \pi  i \,k x) \cos (2 \pi  n x) \, dx &=\frac{i k}{\pi  k^2-\pi  n^2}, \qquad  \forall k\neq \pm n;\\
2 \int_0^1 (1-x) \exp (2 \pi  i n x) \cos (2 \pi  n x) \, dx&=\frac{1}{2}+\frac{i}{4 \pi  n},\qquad \forall n\neq 0.
\end{align*}
Indeed these equalities show that the Fourier coefficient $a_k$ of $\cD(x)$ appears in the diagonal term $q_{n,n}$ as, for $n\neq 0$, 
$$
q_{n,n}=\sum_{k\neq \pm n} a_k\frac{i k}{\pi  k^2-\pi  n^2}+\frac 12(a_n+a_{-n})+\frac{1}{4\pi}(\frac i n a_n-\frac i n a_{-n}).
$$
Since the distribution $\cD(x)$ is real valued, one has $a_{-k}=\overline a_k$ so that the terms $ik a_k$ and $i(-k)a_{-k}=\overline{ik a_k}$ add up to a real contribution for $k\neq \pm n$. With $a_k=x_k+i y_k$ for  $k>0$, one gets for $k\neq n$, 
$$
ika_k+i(-k)a_{-k}=-2ky_k
$$
which gives the first contribution to $q_{n,n}$ as 
$$
\sum_{k>0,k\neq n}y_k\frac{2 k}{\pi  n^2-\pi  k^2}.
$$
 For $k=n$ one gets the two terms
 $$
 \frac 12(a_n+a_{-n})+\frac{1}{4\pi}(\frac i n a_n-\frac i n a_{-n})=x_n-y_n\frac{1}{2\pi n}.
 $$
  For $n=0$ 
$$
q_{0,0}=\sum_{k\neq 0}a_k\frac{i k}{\pi  k^2}+a_0=x_0-\sum_{k> 0}y_k\frac{2}{\pi  k}.
$$ 
\begin{prop}\label{freedom} Let $N\in \N$. The matrices $(q_{i,j})$, $i,j\in \{-N,\ldots,N\}$ obtained from distributions $\cD$ by Proposition \ref{matrixcomp}, are all matrices of the following form, where $a_i$ and $b_j$ are real numbers with $a_{-i}=a_i$ and $b_{-i}=-b_i$, $\forall i\in \{-N,\ldots,N\}$
	\begin{equation}\label{form}
	q_{i,i}=a_i,\quad  \forall i, \qquad  q_{i,j}=\frac{b_i-b_j}{i-j}, \quad \forall j\neq i;
  \qquad i,j\in \{-N,\ldots,N\}.	\end{equation}
  Moreover given a matrix $Q=(q_{i,j})$, $i,j\in \{-N,\ldots,N\}$ of the above form there exists a unique real distribution $\cD(x)$ all of whose Fourier components $a_n=0$ for $n \notin \{-N,\ldots,N\}$ and whose associated matrix is $Q$.
\end{prop}
\proof Proposition \ref{matrixcomp} shows that the matrices $(q_{i,j})$, $i,j\in \{-N,\ldots,N\}$ obtained from distributions $\cD$ by Proposition \ref{matrixcomp} are of the form given by \eqref{form}. Let us show that conversely any matrix of the form \eqref{form} appears.\newline The equality $\psi(x)=\frac 1 \pi \int_0^L \sin(2\pi x(1- y/L))\cD(y)dy$ of Proposition \ref{matrixcomp}, gives, in terms of the Fourier coefficients $a_k=x_k+i y_k$ of $\cD(x)$, for $n\in \Z$, taking $L=1$ for simplicity, 
$$
\psi(n)=-\frac 1 \pi \int_0^1 \sin(2\pi ny)\cD(y)dy=\frac{1}{2\pi i}a_n-\frac{1}{2\pi i}a_{-n}=\frac 1 \pi y_n.
$$ 
Thus the matrix entries  $q_{i,j}$ for $i\neq j$ determine the real numbers $y_j$ which are the imaginary parts of the Fourier coefficients $a_k$ of $\cD$. These imaginary parts then give the following contributions to the diagonal value, for $n\neq 0$,
\begin{equation}\label{diagval}
q_{n,n}=x_n-y_n\frac{1}{2\pi n}+\sum_{k>0,k\neq n}y_k\frac{2 k}{\pi  n^2-\pi  k^2}
\end{equation} 
and we can choose the real part $x_n$ of the Fourier coefficients $a_n$ of $\cD$ to obtain arbitrary diagonal values $a_i$ as required. Moreover the off-diagonal values determine uniquely the imaginary parts $y_n$ of the Fourier coefficients of $\cD(x)$ for $n\in \{-N,\ldots,N\}$ and the diagonal values then determine uniquely the real parts $x_n$ of the Fourier coefficients of $\cD(x)$ for $n\in \{-N,\ldots,N\}$.
\endproof

\begin{rem}\label{nuance} It is important to take as a starting point a distribution $\cD$ on $[0,L]$ and then define the associated quadratic form using \eqref{qua0}, rather than starting from an even distribution on $[-L,L]$. For instance the derivative $\delta'_0$ of the Dirac distribution at $0\in [0,L]$ gives rise to a non-zero quadratic form while the associated even distribution obtained by symmetrisation is equal to $0$.
\end{rem}

\subsection{Relation to the spectral action}
There is a close analogy between the quadratic form in \eqref{form} and the second derivative with respect to perturbations in the spectral action \cite{CC96}, as we will now explain. 
Suppose we are given a linear self-adjoint operator $D$ in a finite-dimensional Hilbert space, which is assumed to have simple spectrum labeled by $\{ \lambda_j\}_{j=-N}^N$, with corresponding eigenvectors $\{e_j\}$.
Suppose that we have a real symmetric positive matrix $Q=(q_{ij})$, $i,j\in \{-N,\ldots,N\}$ defined as
\begin{equation}\label{form2}
  q_{ij} = \left\{ \begin{array}{ll} \frac{b_i - b_j}{\lambda_i -\lambda_j }  & i \neq j \\
a_i & i = j \end{array}
    \right.
\end{equation}
for some $a_i,b_i \in \R$.
The quadratic form $Q$ is given in terms of the Hilbert space inner product as 
\begin{equation}\label{Qquad}
	Q(f,g)=\langle Q f\mid g\rangle=\langle  f\mid Q g\rangle.
\end{equation}

Consider now an even smooth function $f: \R \to \R$. Under perturbations $D \mapsto D+A$ it can be computed that the quadratic form given by the second G\^ateaux derivative of the spectral action is of the form \cite{Han06,Sui11,Skr13,NSkr21,NS21}
$$
\frac 12 \frac{d^2}{dt^ 2} \left( \tr f(D+t A) \right)|_{t=0} = \sum_{i,j} A_{ij} A_{ji} q_{ij}
$$
where $q_{ij}$ is exactly as in \eqref{form} for $b_i = f' (\lambda_i), a_i = f'' (\lambda_i)$ and $\lambda_i = i$, the spectrum of the Dirac operator $D=D_{S^1}$ on the circle.

\section{Finite dimensional even case}

In this section we deal with the general finite dimensional even case.  We are given a real symmetric positive  matrix $Q=(q_{i,j})$, $i,j\in \{-N,\ldots,N\}$ of the form \eqref{form}, {\em i.e.}
\begin{equation*}
	q_{i,i}=a_i,\quad  \forall i, \qquad\ q_{i,j}=\frac{b_i-b_j}{i-j}, \quad \forall j\neq i;
  \qquad i,j\in \{-N,\ldots,N\}	
  \end{equation*}
where the scalars $a_i$ fulfill $a_{-j}=a_j$ and $b_{-j}=-b_j$  for all $j\in \{-N,\ldots,N\}$.

We let $e_j$, $j\in \{-N,\ldots,N\}$, be the canonical orthonormal basis given by the vectors $(\delta(k,j))$ whose all components are $0$ except one. Using the canonical inner product $\langle \alpha\mid \beta\rangle$, the quadratic form $Q$ is given by 
\begin{equation}\label{Qquad}
	Q(f,g)=\langle Q f\mid g\rangle=\langle  f\mid Q g\rangle.
\end{equation}

\begin{lem}\label{basics} $(i)$~Let $\gamma$ be such that $\gamma(e_j):=e_{-j}$ for all $j\in \{-N,\ldots,N\}$. One has $\gamma^2=\id$ and $Q\gamma=\gamma Q$.\newline
$(ii)$~Let $D$ be defined by $D(e_j):=j\,e_{j}$ for all $j\in \{-N,\ldots,N\}$. One has $D\gamma=-\gamma D$ and 
\begin{equation}\label{QD}
D\,	Q-Q\, D=\vert \beta\rangle\langle \eta\vert -\vert \eta\rangle\langle \beta\vert, \ \ \beta=\sum b_j\,e_j,\ \ \eta=\sum e_j.
\end{equation}	
\end{lem}
\proof $(i)$~One has $q_{-i,-j}=q_{i,j}$ for all $i,j\in \{-N,\ldots,N\}$.\newline
$(ii)$~The diagonal elements of the diagonal matrix $D$ are antisymmetric which gives $D\gamma=-\gamma D$. Let us prove \eqref{QD}. One has $(DQ)_{i,j}=iq_{i,j}$, $(QD)_{i,j}=jq_{i,j}$ so that $(D\,	Q-Q\, D)_{i,j}=b_i-b_j$ for all $i,j\in \{-N,\ldots,N\}$. Similarly one has 
$$
(\vert \beta\rangle\langle \eta\vert)_{i,j}=\vert \beta\rangle_i\langle \eta\vert_j=b_i, \ \ (\vert \eta\rangle\langle \beta\vert)_{i,j}=\vert \eta\rangle_i\langle \beta\vert_j=b_j
$$
which gives the required equality.\endproof 
\begin{lem}\label{key} Assume $Q\geq 0$ and $Q\xi=0$ where $\gamma \xi=\xi$ and $\langle \xi\mid \eta\rangle=1$. \newline
$(i)$~One has $Q\, D\,\xi=-\beta $.\newline
$(ii)$~The operator $D':=D-\vert D\,\xi\rangle\langle \eta\vert$ is selfadjoint with respect to the inner product defined by $Q$.	
\end{lem}
\proof $(i)$~We apply \eqref{QD} and get, using  $Q\xi=0$ and $\langle \beta\vert \xi \rangle=0$ since the two eigenspaces of $\gamma$ are orthogonal,
$$
-Q\, D\,\xi=(D\,	Q-Q\, D)\xi=\vert \beta\rangle\langle \eta\vert \xi\rangle -\vert \eta\rangle\langle \beta\vert \xi \rangle=\beta.
$$
$(ii)$~The inner product defined by $Q$ is given by \eqref{Qquad}, \ie 
$$
\langle f\mid g\rangle_Q=\langle Q f\mid g\rangle.
$$
Thus we want to show that $\langle D' f\mid g\rangle_Q=\langle f\mid D'g\rangle_Q$ for all $f,g$. One has, with $R=-\vert D\,\xi\rangle\langle \eta\vert$ 
$$
\langle D' f\mid g\rangle_Q=\langle Q D'f\mid g\rangle=\langle Q Df\mid g\rangle+\langle Q Rf\mid g\rangle.
$$
By $(i)$, one has $QR=-\vert QD\xi\rangle \langle \eta\vert=\vert \beta\rangle\langle \eta\vert$. Moreover by \eqref{QD}, one has $QD-DQ=-\vert \beta\rangle\langle \eta\vert +\vert \eta\rangle\langle \beta\vert$. Thus 
$$
\langle D' f\mid g\rangle_Q=\langle DQf\mid g\rangle+\langle R'f\mid g\rangle, \ \ R'=\vert \eta\rangle\langle \beta\vert.
$$
 Moreover, using that both $Q$ and $D$ are selfadjoint,  
 $$
 \langle f\mid D'g\rangle_Q=\langle Qf\mid D g\rangle+\langle Qf\mid R g\rangle=\langle DQf\mid g\rangle+\langle f\mid QR g\rangle
 $$
 and the required equality follows from 
 $$
 \langle f\mid QR g\rangle=\langle f\mid  \beta\rangle\langle \eta\vert g\rangle=\langle R'f\mid g\rangle.
 $$
 \endproof 
 \begin{lem}\label{lem-main} Let $Q$, $D$, $\xi$, $\eta$ and $D'$ be as in Lemma \ref{key}. Then \newline
 $(i)$~Let $s\notin \{-N,\ldots,N\}$. Then 
 \begin{equation}\label{equiv}
 {\rm Det}(D'-s)=0\iff \sum_{j=-N}^N \frac{\xi_j}{s-j}=0.
\end{equation}
 $(ii)$~One has $ {\rm Det}(D')=0$, and for $j\in \{-N,\ldots,N\}$, $j\neq 0$, ${\rm Det}(D'-j)=0\iff \xi_j=0$. 	
 \end{lem}
 \proof 
 {\em (i)}
 We start by writing, in terms of $R =- \langle D \xi \rangle \langle \eta |$:
 $$
D' -s = D+ R-s = (D-s) \left (\id + (D-s)^{-1} R\right)
$$
Consequently 
$$
\Det(D' - s) = \Det(D-s) \Det (\id + (D-s)^{-1} R).
$$
To compute the second determinant we use the identity 
$$
 \Det (\id + A )=\sum_{k=0}^{\infty} \tr\left(\wedge^k A\right)
$$
applied to the rank one operator $A= (D-s)^{-1} R$. The  higher exterior powers $\wedge^k A$ vanish for $k>1$ thus
$$
 \Det (\id + (D-s)^{-1} R )= 1 - \tr \left( | (D-s)^{-1} D \xi \rangle \langle \eta | \right) = -s \langle \eta | (D-s)^{-1} \xi \rangle, 
 $$
 using $(D-s)^{-1} D \xi =\xi+ s (D-s)^{-1}\xi$ and $ \langle \eta | \xi\rangle=1$.
 Hence,
 \begin{equation}
   \Det(D' - s) = -s \Det (D-s) \langle \eta | (D-s)^{-1}\xi  \rangle = -s \prod_{i=-N}^N (i -s) \sum_{j=-N}^N (j -s)^{-1} \xi_j.
   \label{eq:detDp}
 \end{equation}
 We conclude that if $s \neq j$ for $j = -N, \ldots, N$ then $\Det (D' -s)=0$ iff $\sum_{j=-N}^ N (j -s)^{-1} \xi_j =0$.

 For {\em (ii)} we have that $D ' \xi = 0$ so $\Det (D' ) = 0$.
 For $s= j \neq 0$ we find that the only non-vanishing term on the right-hand side of \eqref{eq:detDp} is $j \prod_{i\neq j} (i -j ) \xi_j$ which is zero iff $\xi_j=0$.
 \endproof 
 \begin{rem} The expression \eqref{eq:detDp} for ${\rm Det}(D'-s)$ is related to 
 the ordinary Lagrange interpolation polynomial for the function $f(\lambda)$ at the points $\lambda_0, \lambda_1, \ldots, \lambda_n$   :
$$
P(x) = \sum_{k=0}^{n} f(\lambda_k) \cdot \prod_{\substack{j=0 \\ j \ne k}}^{n} \frac{x - \lambda_j}{\lambda_k - \lambda_j}.
$$
 \end{rem}

We now compute the Fourier transform of functions on $[0,L]$ translated to $[-\frac L2,\frac L2]$ and extended by $0$ to the full line $\R$. The Fourier transform is defined by
$$
\fourier(f)(s):=\int_\R f(x)\exp(- isx)dx.
$$
The next Proposition is a reformulation of the Shannon sampling theorem (\cite{Shannon}).
\begin{prop}Let $f\in L^2([0,L])$ and $f^\sigma(x):=f(x+\frac L2)$ for $\vert x\vert\leq \frac L2$ be extended by $0$ on $\R$.\newline
$(i)$~The restriction of the Fourier transform of $f^\sigma$ to $\frac{2\pi}{L}\Z$ is given by the Fourier transform $\widehat{ f}$ of $f\in L^2(\R/L\Z)$ as follows
\begin{equation}\label{fourz}
\fourier(f^\sigma)(\frac{2\pi}{L}n)=(-1)^n\, \widehat{ f}(n), \qquad \forall n\in \Z.
\end{equation} 
$(ii)$~The Fourier transform of $f^\sigma$ is given by 
\begin{equation}\label{four}
\fourier(f^\sigma)(s)=\sin(L s/2)\sum_\Z  \widehat{ f}(n)\frac{1}{Ls/2-n\,\pi }.
\end{equation}	
\end{prop}
\proof $(i)$~Let $n\in \Z$. One has by definition $\widehat{ f}(n)=\int_{0}^{L}f(x)\exp(-2\pi inx/L)dx$ and
$$
\fourier(f^\sigma)(\frac{2\pi}{L}n)=\int_\R f^\sigma(x)\exp(-2\pi inx/L)dx=\int_{-\frac L2}^{\frac L2}f(x+\frac L2)\exp(-2\pi inx/L)dx=$$ $$(-1)^n\,\int_{0}^{L}f(x)\exp(-2\pi inx/L)dx=(-1)^n\,\widehat{ f}(n).
$$
$(ii)$~One has $f(x)=\frac 1L\sum_\Z \widehat{ f}(n)\exp(2\pi inx/L)$, thus it is enough to treat the case $f(x)=\exp(2\pi inx/L)$. Then 
$$
\fourier(f^\sigma)(s)=\int_\R f^\sigma(x)\exp(- isx)dx=\int_{-\frac L2}^{\frac L2}\exp(2\pi in(x+\frac L2)/L)\exp(- isx)dx=$$ $$=(-1)^n\,\int_{-\frac L2}^{\frac L2}\exp(i(\frac{2\pi}{L}n-s)x)dx=(-1)^n\,\frac{1}{i(\frac{2\pi}{L}n-s)}(\exp(i(\frac{2\pi}{L}n-s)x))_{-\frac L2}^{\frac L2}=2L\frac{\sin(Ls/2)}{Ls-2\pi n}
$$
which gives \eqref{four} using $f(x)=\frac 1L\sum_\Z \widehat{ f}(n)\exp(2\pi inx/L)$. \endproof 
After these preliminaries we obtain
 \begin{thm}\label{finmain} Let $Q$ be a real symmetric positive  matrix  of the form \eqref{form} with one dimensional kernel which is even with respect to $\gamma$. Let $\xi\in \ker\, Q$.\newline
 $(i)$~All the roots of the following polynomial are real:
\begin{equation}\label{thmpoly}
 P(s)=\sum_{k\in\{-N,\ldots,N\}}\xi_k\times \left(\prod_{j\in\{-N,\ldots,N\}, j\neq k}(j-s)\right).
\end{equation} 	
$(ii)$~The Fourier transform $\widehat\xi(z)$ of the function 
  $$
\xi(x):=\sum \xi_k\exp(2\pi ikx), \ \forall x\in [0,1], \ \xi(x)=0, \ \forall x\notin [0,1] 
$$
is entire and has all its zeros on the real line.
 \end{thm}
 \proof {\em (i)}
Recall $D, \gamma$ from Lemma \ref{basics} and note that $D$ has one-dimensional kernel. If  $e_0\in \ker\, Q$ one checks $(i)$ and $(ii)$ directly. Thus we assume that $e_0\notin \ker\, Q$. Let $\xi\in \ker\, Q$, $\xi\neq 0$. The real symmetric positive  matrix $Q$ defines an inner product on $\R^{2N+1}$ and its radical consists of the one dimensional subspace generated by $\xi$. One has $D\xi\neq 0$ since otherwise $e_0\in \ker\, Q$.  One has $QD\xi\neq 0$ since $D\xi$ is odd and therefore linearly independent of $\xi$ while $\ker\, Q$ is one-dimensional. By \eqref{QD} one has
 $$
 0\neq (D\,	Q-Q\, D)(\xi)=\vert \beta\rangle\langle \eta\vert\xi\rangle -\vert \eta\rangle\langle \beta\vert\xi\rangle=\vert \beta\rangle\langle \eta\vert\xi\rangle.
 $$ 
 Thus one can normalize $\xi$ so that $\langle \eta\vert\xi\rangle=1$. Let then $D':=D-\vert D\,\xi\rangle\langle \eta\vert$ as in Lemma \ref{key}. One has $D'(\xi)=0$ and thus $D'$ induces an operator $D"$ on the Euclidean space $E$ obtained as the separated quotient of $(\R^{2N+1},Q)$. By Lemma \ref{key}, $(ii)$, the operator $D"$ is selfadjoint in $E$. Thus the real spectral theorem (see \cite{AS}, Theorem 7.29) applies and shows that the characteristic polynomial of $D"$ has all its roots in $\R$. Let $v_j$ be an orthonormal basis of $E$ of eigenvectors for  $D"$ with eigenvalues $\lambda_j$. Let $w_j\in \R^{2N+1}$ be lifts of the $v_j$. One has $D"(v_j)=\lambda_jv_j$ and hence $D'(w_j)=\lambda_jw_j+s_j \xi$ for some real scalars $s_j$. Thus in the basis of $\R^{2N+1}$ formed by $\xi$ and the $w_j$, the matrix of $D'$ is triangular, with $0$ and the $\lambda_j$ on the diagonal. Thus  one gets 
 $$
 {\rm Det}(D'-s)=-s\prod (\lambda_j-s).
 $$
	Comparing this formula with \eqref{eq:detDp} one obtains that the polynomial $P(s)$ of \eqref{thmpoly} has all its zeros in $\R$.\newline
	$(ii)$~The Fourier transform of the function with support in $[0,1]$ given there by $\exp(2\pi ikx)$ is 
	$$
	\int_0^1 \exp (2 \pi  i k x) \exp (-i s x) \, dx=\frac{2 e^{-\frac{is}{2} } \sin \left(\frac{s}{2}\right)}{s-2 \pi  k}.
	$$
	Thus the  Fourier transform of $\xi(x)$ is 
	$$
	\widehat\xi(z)=2 e^{-\frac{iz}{2} } \sin \left(\frac{z}{2}\right)\left(\sum_{\{-N,\ldots,N\}} \frac{\xi_j}{z-2\pi j}\right).
	$$
	The zeros $z\in 2\pi \Z$ of $\sin \left(\frac{z}{2}\right)$ cancel the pole at $2\pi j$ which occurs when $\xi_j\neq 0$ and remain as zeros of $\widehat\xi(z)$ otherwise. The other zeros are given by the roots of $P(z/2 \pi)=0$ where $P(z)$ is defined in \eqref{thmpoly}. Thus by $(i)$, all these zeros are real. \endproof 
	\endproof

\subsection{General finite dimensional operator $D$}
Consider as before the more general situation of a linear self-adjoint operator $D$ in a finite-dimensional Hilbert space, which is assumed to have simple spectrum labeled by $\{ \lambda_j\}_{j=-N}^N$, with corresponding eigenvectors $\{e_j\}$. We then consider the real symmetric positive matrix $Q=(q_{ij})$, $i,j\in \{-N,\ldots,N\}$ defined in \eqref{form2}. 

The following Lemma's are the analogues of Lemmas \ref{basics}, \ref{key} and \ref{lem-main}, whose proofs follow {\em mutatis mutandis}.
\begin{lem}
  \label{basics-general}
  Suppose that $\lambda_{-i} = - \lambda_i$ and $b_{-i} = -b_{i}$ for all $i \in \{ -N, \ldots, N\}$.\\
$(i)$~Let $\gamma$ be such that $\gamma(e_i):=e_{-i}$ for all $i\in \{-N,\ldots,N\}$. One has $\gamma^2=\id$ and $Q\gamma=\gamma Q$.\newline
$(ii)$ One has $D\gamma=-\gamma D$ and 
\begin{equation}\label{QD-sa}
D\,	Q-Q\, D=\vert \beta\rangle\langle \eta\vert -\vert \eta\rangle\langle \beta\vert, \ \ \beta=\sum b_i \,e_i,\ \ \eta=\sum e_i,
\end{equation}
so that $\beta$ is odd and and $\eta$ is even with respect to the $\Z_2$-grading given by $\gamma$.
\end{lem}

\begin{lem}
  \label{key-general}
Let $D,Q,\gamma$ be as in Lemma \ref{basics-general}, assume $Q\geq 0$ and $Q\xi=0$ where $\gamma \xi=\xi$ and $\langle \xi\mid \eta\rangle=1$. \newline
$(i)$~One has $Q\, D\,\xi=-\beta $.\newline
$(ii)$~The operator $D':=D-\vert D\,\xi\rangle\langle \eta\vert$ is self-adjoint with respect to the inner product defined by $Q$.	
\end{lem}

 \begin{lem}\label{main-odd} Let $Q$, $D$, $\xi$, $\eta$ and $D'$ be as in Lemmas \ref{basics-general} and \ref{key-general}, and assume that $D$ has simple spectrum. Then
   \begin{enumerate}[(i)]
\item Let $s\in \C \setminus \{-\lambda_N,\ldots,\lambda_N\}$ and $s \neq 0$. Then 
 \begin{equation}\label{equiv}
 {\rm Det}(D'-s)=0\iff \sum_{j=-N}^N \frac{\xi_j}{s-\lambda_j}=0.
\end{equation}
\item One has $ {\rm Det}(D')=0$, and if $\lambda_j \neq 0$ we have ${\rm Det}(D'-\lambda_j )=0\iff \xi_j=0$.
   \end{enumerate}
   \end{lem}

As a result, we also have the following analogue of Theorem \ref{finmain}(i):
\begin{prop}\label{prop:finmain} Let $D$ have simple spectrum and let $Q$ be a real symmetric positive  matrix  of the form \eqref{form2} with one dimensional even kernel. Let $\xi\in \ker\, Q$. Then all the roots of the following polynomial are real:
\begin{equation}\label{proppoly}
 P(s)=\sum_{k\in\{-N,\ldots,N\}}\xi_k\times \left(\prod_{j\in\{-N,\ldots,N\}, j\neq k}(\lambda_j-s)\right).
\end{equation} 	
 \end{prop}
 We then obtain Theorem \ref{finmain}(ii) as a corollary to this result, when it is applied to the case $D= D_{S^ 1}$ so that $\lambda_j = j$.

\section{Infinite dimensional case}
In this section we prove the result announced in the Introduction, 
\begin{thm}\label{main}
Let $ L > 0 $, and let $\cD$ be a real  distribution on the interval $ [0, L] $. Let $Q$ be the quadratic form defined on trigonometric polynomials by \eqref{qua0}. Assume that $Q$ defines  a lower-bounded essentially selfadjoint operator   and that the minimum of its spectrum is a simple, isolated eigenvalue \( \lambda \), with even eigenfunction\footnote{\ie invariant under the symmetry $x\mapsto L-x$ of $[0,L]$} $\xi$. Then all the zeros of the entire function $\widehat \xi(z)$, $z\in \C$, the Fourier transform  of $\xi$ lie on the real line.
\end{thm}
\proof By hypothesis the trigonometric polynomials form a core for the operator $A$ which is defined by
$$
\langle \alpha\mid A \beta\rangle=\langle \alpha\mid \beta\rangle_Q.
$$
 We normalize the eigenvector $\xi$ by $\Vert \xi\Vert=1$.   Given $\epsilon>0$ there exists an even trigonometric polynomial $\eta_\epsilon$ with 
 \begin{equation}\label{approx}
 	\Vert \eta_\epsilon\Vert=1, \ \ \Vert\xi-\eta_\epsilon\Vert<\epsilon, \  \ \langle \eta_\epsilon\mid \eta_\epsilon\rangle_Q<\lambda+\epsilon.
 \end{equation}
 Let $N$ be a finite integer such that the support of $\eta_\epsilon$ is contained in $\{-N,\ldots,N\}$. By Proposition \ref{matrixcomp}  the matrix of the restriction $Q_N$ of the quadratic form $Q$ to the space $E_N$ of trigonometric polynomials with support  in $\{-N,\ldots,N\}$ is of the form \eqref{form}. Since $Q_N$ is a restriction of $Q$ to a subspace, its minimum is $\geq \lambda$ and the above inequalities show that it is between $\lambda$ and $\lambda+\epsilon$. By hypothesis the spectrum of the operator $A$ is contained, except for the simple eigenvalue $\lambda$ in the interval $[\lambda+\delta, \infty)$ for some $\delta>0$. Thus the restriction of $Q$ to the orthogonal complement of $\xi$ fulfills
 $$
 \langle \alpha\mid \alpha\rangle_Q\geq (\lambda+\delta)\Vert \alpha\Vert^2, \ \ \forall  \alpha \mid \langle \alpha\mid \xi\rangle=0.
 $$
 This holds in particular in the codimension one subspace $F_N$ of $E_N$ defined by $\langle \alpha\mid \xi\rangle=0$ (note that $F_N$ cannot be all of $E_N$ since then we also would have had $\xi \perp \eta_\epsilon$ which contradicts \eqref{approx}). The smallest eigenvalue of $Q_N$ fulfills $\lambda\leq \lambda_N\leq \lambda+\epsilon$. By the  min-max theorem, the next eigenvalue $\mu_N\geq \lambda_N$  of $Q_N$ is given by 
$$
\mu_N=\max _{\substack{\mathcal{M} \subset E_N \\ \operatorname{codim}(\mathcal{M})=1}} \min _{x \in \mathcal{M}, \Vert x\Vert=1}Q_N(x)
$$
so that, using $\mathcal{M}=F_N$ one gets $\mu_N\geq \lambda+\delta$. This shows that for $\epsilon<\delta/2$ the eigenvalue $\lambda_N$ of $Q_N$ is simple and the only one in the interval $[\lambda,\lambda+\delta]$. Let then $P$ be the spectral projection of $Q_N$ for the eigenspace associated to the minimal eigenvalue $\lambda_N$. Decomposing 
$$\eta_\epsilon= P\eta_\epsilon+ (1-P)\eta_\epsilon=\alpha+\beta \ \Rightarrow Q(\eta_\epsilon)=\lambda_N \Vert\alpha\Vert^2+Q(\beta)$$ 
where $Q(\beta)\geq (\lambda+\delta)\Vert\beta\Vert^2$. Thus 
 by \eqref{approx}, one gets that the convex combination with weights $\Vert\alpha\Vert^2$ and $\Vert\beta\Vert^2$ of $\lambda$ and $Q(\beta)/\Vert\beta\Vert^2\geq (\lambda+\delta)$ is less than $\lambda+\epsilon$. It follows that $\Vert\beta\Vert^2=1-\Vert\alpha\Vert^2\leq \epsilon/\delta$. Let then $\xi_N$  be the eigenvector of $Q_N$ for the eigenvalue $\lambda_N$ given by $P\eta_\epsilon$. One controls 
 $$
 \Vert\eta_\epsilon-P\eta_\epsilon\Vert\leq \sqrt{\epsilon/\delta}, \  \ \Vert\xi-\eta_\epsilon\Vert<\epsilon\Rightarrow \Vert \xi -\xi_N\Vert \leq \epsilon+\sqrt{\epsilon/\delta}.
 $$
 Since the even functions form a closed subspace in $L^2[0,L]$, and $\xi$ is even, for $\epsilon$ small enough the vector $\xi_N$ is even as well. Thus, it follows from Theorem \ref{finmain} that all the zeros of the Fourier transform $\widehat \xi_N$ are real. Moreover when $\epsilon\to 0$ the vectors $\xi_N$ converge in norm to $\xi$ so that the sequence $\widehat \xi_N(z)$ converges to $\widehat \xi(z)$ uniformly on compact subsets in $\C$ (as in the proof of Theorem \ref{thmtoepcont}). But then the Hurwitz theorem applies, allowing us to conclude that all zeros of $\widehat \xi(z)$ are real.
 \endproof

 \section{Spectral action and divided differences}
 \label{sect:sa}

 As already observed there is a close relation of the quadratic form in \eqref{Qquad} and the spectral action $\tr f(D)$ introduced in \cite{CC96}. We will now analyze this in more detail for perturbations of the type  $D \mapsto D +R$ with $R = - | D \xi \rangle \langle \eta|$ as in Lemma \ref{key}, extending the Taylor expansions derived in \cite{Han06,Sui11,Skr13,NSkr21,NS21} to this case.

First, in order to make sense of the spectral action for (not necessarily self-adjoint or normal) bounded perturbations of a self-adjoint operator we write $f (x)$ as a Fourier transform, and then invoke Araki's expansionals \cite{araki} ---or Dyson series--- to give meaning to $e^{i (H_0+ V)}$ for bounded perturbation $V$ of a self-adjoint operator $H_0$ (to be precise, this is \cite[Eq. 5.16]{araki}):
\begin{equation}
e^{i (H_0 +V)}:=
\mathrm{Exp}_r \left(\int_0^1 ; i  e^{i s H_0}   V e^{-i s H_0} ds \right) e^{i H_0} := \sum_{n \geq 0} i^n \int_{\Delta_n} e^{i s_0 H_0} V e^{i s_1 H_0} \cdots V e^{i s_n H_0}  d^n s.
\label{eq:expans}
\end{equation}
where the $n$-simplex $\Delta_n$ is parametrized by tuples $(s_0,\ldots, s_n) \in \R^{n+1}_{\geq 0}$ satisfying $\sum_k s_k = 1$. In our case of interest, these expansionals are given by series of the following form:
\begin{equation}
e^{i \xi (D +t R)}
:= \sum_{n \geq 0} (it \xi)^n \int_{\Delta_n} e^{i s_0\xi D} R e^{i s_1 \xi D} \cdots R e^{i s_n \xi D}  d^n s.
\label{eq:expans}
\end{equation}
Note that since $|\Delta_n | = 1/n!$ the $n$'th summand in the expansional is norm bounded by $t^n|\xi|^n  \| R\|^n/n!$ so that the series expansion is norm convergent. This suggest to define for suitable functions $f$: 
\begin{equation}
  f( D+ t R) := \int_\R \widehat f(\xi) e^{i \xi (D+t R)} d\xi 
.
  \label{eq:f-R}
\end{equation}
More precisely, we have 
\begin{lem}
  Let $D$ be a self-adjoint operator on $\cH$, $R$ bounded operator on $\cH$ and let $f$ be such that $\| \widehat{f^{(n)}}\|_1 \leq C^{n+1} n!$ for all $n \geq 0$ and some $C\geq 1$. Then for sufficiently small $t$ the expression in \eqref{eq:f-R} is a bounded operator on $\cH$.
\end{lem}
\proof
We estimate:
$$
\|  f( D+ t R) \| \leq \sum_{n \geq 0} \frac{t^n \| R\| ^n}{n!} \int_\R | \widehat f(\xi)  \xi^n | d \xi \leq \frac{C}{1- t C \| R\|}
$$
which is bounded for $t < 1/ (C \| R\|)$.
\endproof

In order to define the spectral action as the trace of this operator, we need a more restrictive class of functions. Namely, as in \cite{NS21} we define
$$
\mathcal E^{s} :=\left\{f\in C^\infty :~ \text{there exists $C\geq 1$ s.t. } \|\widehat{(fu^m)^{(n)}}\|_1\leq C^{n+1}n! \,\text{ for all $m\leq s$ and $n\geq 0$}\right\},
$$
where $u(x) = x-i$.
\begin{lem} If $D$ is $s$-summable, {\em i.e.} $(D-i)^{-s}$ is trace-class for some $s\geq 0$, and $f \in \mathcal E^{s}$ then $f(D+tR)$ is trace-class for sufficiently small $t$.  
  \end{lem}
\proof
The proof of the required estimates follows line-by-line the proof of \cite[Theorem 6]{NS21} after having given the meaning \eqref{eq:expans} to the exponentials appearing in the multiple operator integrals ({\em i.e.} Definition 2 in {\em loc. cit.}). 
\endproof
We recall the definition of divided differences. Let $f: \R \to \R$ and let $x_0, x_1, \ldots x_n$ be distinct points in $\R$. The {\rm divided difference of order $n$} is defined by the recursive relations
\begin{align*}
f[x_0] &= f(x_0), \\
f[x_0,x_1, \ldots x_n] &= \frac{ f[x_1, \ldots x_n] -f[x_0,x_1, \ldots x_{n-1}]}{x_{n} - x_0}.
\end{align*}
Also note the following useful representation, due to Hermite \cite{hermite}: for any $x_0, \ldots, x_n \in \R$,
$$
f[x_0, x_1, \ldots, x_n] = \int_{\Delta_n} f^{(n)} \left(s_0 x_0 + s_1 x_1 + \cdots + s_n x_n\right) d^ns.
$$
This also allows to extend the definition of divided difference to coinciding points.

 \begin{lem}
   \label{lem:sa}
Let $D$ be a self-adjoint operator in $\cH$ such that $(D-i)^{-s}$ is trace class for some $s \geq 0$, $R$ is a bounded operator in $\cH$ and $f \in \mathcal E^{s}$. Then $t \mapsto \tr f(D+ tR)$ is smooth in a neighborhood of $0$ with $n$'th derivative
   $$
   \frac{d^n}{dt^n} \tr f(D+tR) | _{t=0} = n! \sum R_{i_1i_2} \cdots R_{i_n i_1} f' [\lambda_{i_1},\ldots, \lambda_{i_n} ]. 
   $$
   in terms of the eigenvalues $\lambda_i$ of $D$.
   \end{lem}
 \proof
Since $f(D+tR)$ is defined in terms of Araki's expansional formula, we have 
\begin{align*}
  \frac{d^n}{dt^n} \tr f(D+tR) |_{t=0} &=   \int \frac{d^n}{dt^n} \tr ( \widehat f (\xi) e^{i \xi (D+tR)} )|_{t=0}  d \xi  \\
  &= n! \sum_{i_1,i_2,\ldots, i_n} R_{i_1 i_2} \cdots R_{i_n i_1} \int (i \xi)^n \exp[i\xi \lambda_{i_1},\ldots ,i\xi \lambda_{i_n},i\xi \lambda_{i_1}] \widehat f(\xi) d\xi\\
  &= n! \sum_{i_1,i_2,\ldots, i_n} R_{i_1 i_2} \cdots R_{i_n i_1} f[ \lambda_{i_1},\ldots , \lambda_{i_n},\lambda_{i_1}].\qedhere
  \end{align*} 
\endproof

\begin{rem}
  It would be interesting to extend the above definition of the spectral action for not necessarily self-adjoint perturbations to the case where also the assumption on self-adjointness on the operator $D$ is relaxed. This has potential applications to Lorentzian spectral triples, 
  \end{rem}

 Let us now take a function $f$ such that $f''(\lambda_j) = a_j$ and $f' (\lambda_j) = b_j$, where $a_j,b_j$ are the coefficients of the quadratic form as in \ref{form2}. We then have 
$$
 \frac{d}{dt} \tr f(D+tR) |_{t=0} = \sum R_{jj} b_{j}; \qquad 
 \frac{d^2}{dt^2} \tr f(D+tR) |_{t=0} = \sum R_{ij} R_{ji} q_{ij}. 
$$

 \begin{prop}
   \label{prop:der-sa}
   Let $\cH$ be finite-dimensional. In the notation of Lemma \ref{basics-general}, assume $Q= (q_{ij}) \geq 0$ and $Q \xi =0$ where $\gamma \xi =\xi$ and $\langle \xi | \eta \rangle = 1$. 
   Let $R = - | D \xi \rangle \langle \eta |$ so that $R_{ij} = - (D\xi)_i$. Then we have 
   $$
 \frac{d}{dt} \tr f(D+tR) |_{t=0} = \langle D \xi, D \xi \rangle_Q; \qquad 
 \frac{d^2}{dt^2} \tr f(D+tR) |_{t=0} =  \langle D \xi, D \xi \rangle_Q.
$$
 \end{prop}
 \proof
 The second derivative is reduced to the first derivative because  $QD \xi = -\beta$ ({\em cf.} Lemma \ref{key-general}(i)). Indeed, this yields:
 $$
\sum_j R_{ji} f' [\lambda_{i},\lambda_{j}] = - \sum_j f' [\lambda_{i},\lambda_{j}] (D \xi)_j = - (Q D \xi)_i  = b_i = f'(\lambda_i). 
$$
From this it follows that
\[
\sum R_{ij} R_{ji} f' [\lambda_{i},\lambda_{j}] = \sum_i  R_{ii} f'(\lambda_i) =-  \sum_i  (D \xi)_i f'(\lambda_i) =  \sum_i  (D \xi)_i \overline{ (Q D \xi)_i} = \langle D \xi , D \xi \rangle_Q.
\qedhere
\]
\endproof


\appendix

\section{An instance of truncation matrices}\label{double}
In this appendix we  describe an example  where the truncation matrices admit simple maximal and minimal eigenvalues but this property fails in the limit where the maximal and minimal eigenvalues have multiplicity $2$. 
We let $L=1$ and take the distribution $\cD$  of the form 
\begin{equation}\label{simple}
\cD(x)= \delta_0(x)+2\pi\, b\, \sin (2 \pi  x).
\end{equation}
We then compute $\psi(x):=\frac 1 \pi \int_0^1 \sin(2\pi x(1- y))\cD(y)dy$ as 
$$
\psi(x)=\frac{1}{\pi }\sin (2 \pi  x)+2b \int_0^1 \sin (2 \pi  y) \sin (2 \pi  x (1-y)) \, dy=\frac{1}{\pi }\sin (2 \pi  x)+b\frac{  \sin (2\pi  x)}{\pi -\pi  x^2}=\frac{  \sin (2\pi  x)}{\pi}\left(1+\frac{b}{1-x^2}\right).
$$
Thus one has $\psi(n)=0$ for all $n\in \Z$, except for $n=\pm 1$. Moreover one gets $\psi(-1)=b$ and $\psi(1)=-b$. The derivative is 
$$
\psi'(x)=\frac{2 b\, x\, \sin (2 \pi  x)}{\pi  \left(1-x^2\right)^2}+\frac{2 b\, \cos (2 \pi  x)}{1-x^2}+2\, \cos (2 \pi  x).
$$
One finds that for $n\in \Z$, not equal to $\pm 1$, one has $$\psi'(n)=2 + \frac{2 b\, }{1-n^2}$$
while for $n=\pm 1$ one has $\psi'(n)=\frac{b}{2}+2$.
\begin{lem}\label{matrixmu} Let $L=1$ and $\cD$ be given by \eqref{simple}, $Q$ be the quadratic form $Q$ of \eqref{matelem}.  The matrix elements $q_{n,m}$ of $Q$ are given by $q_{n,m}=2\delta_{n,m}+b\,\mu_{n,m}$, where the matrix $\mu$ is independent of $b$ and given by $\mu_{n,m}=0,\  \forall n,m  \notin \{-1,0,1\},\ n\neq m$ and
$$
(\mu_{n,m})_{n,m\in \{-1,0,1\}}=\left(
\begin{array}{ccc}
 \frac{1}{2} & -1 & -1 \\
 -1 & 2 & -1 \\
 -1 & -1 & \frac{1}{2} \\
\end{array}
\right), \ \ \mu_{n,n}=\frac{2 \, }{1-n^2}, \ \forall n\notin \{-1,0,1\},
$$
$$
\mu_{n,-1}=\mu_{-1,n}=\frac{1}{-1-n}, \ \forall n\notin \{-1,0,1\}, \ \ \mu_{n,1}=\mu_{1,n}=\frac{1}{-1+n}, \ \forall n\notin \{-1,0,1\}.
$$
\end{lem}
\proof This follows from Proposition \ref{matrixcomp} and the above determination of $\psi(n)$ and $\psi'(n)$.\endproof 
The matrix elements $\mu_{n,m}$ of  $\mu$, for $\vert n\vert \leq 4$ and $\vert m\vert \leq 4$ are the following

$$
\left(
\begin{array}{ccccccccc}
 -\frac{2}{15} & 0 & 0 & \frac{1}{3} & 0 & -\frac{1}{5} & 0 & 0 & 0 \\
 0 & -\frac{1}{4} & 0 & \frac{1}{2} & 0 & -\frac{1}{4} & 0 & 0 & 0 \\
 0 & 0 & -\frac{2}{3} & 1 & 0 & -\frac{1}{3} & 0 & 0 & 0 \\
 \frac{1}{3} & \frac{1}{2} & 1 & \frac{1}{2} & -1 & -1 & -\frac{1}{3} & -\frac{1}{4} & -\frac{1}{5} \\
 0 & 0 & 0 & -1 & 2 & -1 & 0 & 0 & 0 \\
 -\frac{1}{5} & -\frac{1}{4} & -\frac{1}{3} & -1 & -1 & \frac{1}{2} & 1 & \frac{1}{2} & \frac{1}{3} \\
 0 & 0 & 0 & -\frac{1}{3} & 0 & 1 & -\frac{2}{3} & 0 & 0 \\
 0 & 0 & 0 & -\frac{1}{4} & 0 & \frac{1}{2} & 0 & -\frac{1}{4} & 0 \\
 0 & 0 & 0 & -\frac{1}{5} & 0 & \frac{1}{3} & 0 & 0 & -\frac{2}{15} \\
\end{array}
\right).
$$

\medskip
Let then $\xi$ be the vector with coordinates $\xi(n)=\frac{1}{1+n}$ for $n\neq -1$ while  $\xi(-1)=0$, and $\eta$ with $\eta(n)=\frac{1}{-1+n}$ for $n\neq 1$, while $\eta(1)=0$. Also let the $U_n$ form the canonical orthonormal basis. Thus one has
$$
\xi=\sum_{n\neq -1}\frac{1}{1+n} U_n, \ \ \eta=\sum_{n\neq 1}\frac{1}{-1+n} U_n.
$$
The functions corresponding to these vectors are 
$$
\xi(x)=-2 \pi  i \left(x-\frac{1}{2}\right) \exp (-2 \pi  i x), \ \ \eta(x)=-2 \pi  i \left(x-\frac{1}{2}\right) \exp (2 \pi  i x).
$$
One has $\overline{\xi(x)}=-\eta(x)$
$$
\sum_{-k}^k\xi(n)\eta(n)=\frac{k^2-3 k-2}{2 k (k+1)}, \  \ \langle \xi\mid \eta\rangle=\frac 12.
$$
We consider the rank $4$ matrix given by 
$$
R:=\vert \eta\rangle \langle U_{1}\vert +\vert U_{1}\rangle \langle\eta \vert-\vert \xi\rangle \langle U_{-1}\vert -\vert U_{-1}\rangle \langle\xi \vert.
$$
Its matrix elements $R_{n,m}$ fulfill $R_{n,m}=0,\  \forall n,m  \notin \{-1,1\}$, while 
$$
R_{n,1}=R_{1,n}=\eta(n), \ \ R_{n,-1}=R_{-1,n}=-\xi(n), \ \ \forall n\notin \{-1,1\}.
$$
The corresponding Schwartz kernel $r(x,y)$ is given by 
$$
r(x,y)=\eta(x)\exp(-2 \pi i y)-\exp(2 \pi i x)\xi(y)-\xi(x)\exp(2 \pi i y)+\exp(-2 \pi i x)\eta(y).
$$
\begin{lem} $(i)$~Let $D$ be the diagonal matrix with diagonal elements $d_n=\frac{2}{1-n^2}$ for $n^2\neq 1$ and $d_n=\frac 12$ for $n^2=1$. One then has $\mu=D+R$.\newline
$(ii)$~The operator $D$ is given by the convolution among periodic functions with period $1$ by the function $$\alpha(x):=-4 \pi (x-\frac 12) \sin(2 \pi x), \ \ \forall x\in [0,1).$$
\end{lem}
\proof $(i)$~This follows from Lemma \ref{matrixmu}.\newline
$(ii)$~One has, for $n\in \Z$, $n\neq \pm 1$, 
$$
2 \pi  \int_0^1 \left(x-\frac{1}{2}\right) \sin (2 \pi  x) \cos (2 \pi  n x) \, dx=\frac{1}{n^2-1}.
$$
while the value of this integral is $-\frac{1}{4}$ for $n=\pm 1$.
 \endproof
 Note that the function $\alpha(x)$ needs to be viewed as a periodic function of period $1$ and this requires reinterpreting the term $x-\frac 12$ as $x-E(x)-\frac 12$ where $E(x)$ is the integer part of $x\in \R$. This plays a role in the formula for the convolution, which is 
 $$
 D(f)(x)=\int_0^1 \alpha(x-y)f(y)dy.
 $$
 One has $$\sin (2 \pi  (x-y))=-\cos (2 \pi  x) \sin (2 \pi  y)+\sin (2 \pi  x) \cos (2 \pi  y)$$
 so that 
 $$
 \alpha(x-y)=-4 \pi (x-y-E(x-y)-\frac 12)\left(-\cos (2 \pi  x) \sin (2 \pi  y)+\sin (2 \pi  x) \cos (2 \pi  y)\right)
 $$
 and the only term which does not separate as a product of functions of $x$ by a function of $y$ is the term involving $E(x-y)$ which is equal to $0$ if $y\leq x$ and to $-1$ if $y>x$. In fact it is more symmetric to add the $\frac 12$ to $E(x-y)$ which coincides with $\frac 12 {\rm Sign}(x-y)$. The other contributions are given by the $4$ terms which are, up to the overall factor $-2 \pi$, 
 $$
 (-i) x e^{-2 i \pi  (y-x)}+i x e^{2 i \pi  (y-x)}-i y e^{2 i \pi  (y-x)}+i y e^{-2 i \pi  (y-x)}
 $$
 which, taking into account the factor $-2 \pi$ can be rewritten as 
 $$
 -\eta(x)\exp(-2i \pi y)+\xi(x)\exp(2i \pi y)-\exp(-2i \pi x)\eta(y)+\exp(2i \pi x)\xi(y)
 $$
 We thus see that these terms cancell the rank $4$ additional contribution $R$ and thus we get the following simple formula for the Schwartz kernel $\mu(x,y)$ of the operator $\mu$,
 \begin{prop} $(i)$~The operator $\mu$ is given by the formula
 $$
 \mu(f)(x)=2 \pi\int_0^1 {\rm Sign}(x-y)\sin(2\pi(x-y))f(y)dy.
 $$
 $(ii)$~In general  the Schwartz kernel of the operator associated to the distribution $\cD$ is equal to the restriction to $x,y\in [0,L]$ of $\cD(\vert x-y\vert)$. 	
 \end{prop}
 \proof $(i)$~Follows from the above computation.\newline
 $(ii)$~By \eqref{quad} one has for smooth $f,g$ with support in $[0,L]$,
 $$
 Q(f,g)=\int_0^L \left((g^**f)(y)+(g^**f)(-y)\right)\cD(y)dy=\int_{-L}^L (g^**f)(y)\cD(\vert y\vert)dy $$
 where one needs to be careful in doubling the coefficient of $\delta_0$ in $\cD(\vert y\vert)$. This formula does not change if one replaces $f(x)$ by $f(x+\frac 12)$ and $g(x)$ by $g(x+\frac 12)$, shifting their supports to the symmetric interval $[-\frac L2,\frac L2]$. Then $Q(f,g)=0$ when $f$ and $g$ have opposite  parity.   One then uses the formula for three  functions $f,g,h$, $h$ even
 $$
 \int_{-L}^L (g^**f)(y)h(-y)dy=\int_{x+y+z=0}g^*(x)f(y)h(z)\omega =\int_{x=y+z}\overline{g(x)}f(y)h(z)dx dy$$
 $$ =\int\overline{g(x)}\int k(x,y)f(y) dy dx,\ \  \ \ k(x,y)=h(x-y)
 $$
where $\omega$ is the measure on the plane $x+y+z=0$ given by $\vert dx\wedge dz\vert=\vert dx\wedge dy\vert$.\endproof 
One can then investigate numerically the zeros of the polynomials  $P_n^\pm(s)$ associated by \eqref{thmpoly} to  the eigenvectors for the largest and smallest eigenvalues of the matrix $\mu(n)$ of size $2n+1$,  which is the compression of the matrix  $\mu$, 
$$
\mu(n)_{i,j}:=\mu_{i,j}, \ \ \forall i,j \in \{-n,\ldots , n\}.
$$
For $n=1$, this matrix is simply 
$$
\mu(1)=\left(
\begin{array}{ccc}
 \frac{1}{2} & -1 & -1 \\
 -1 & 2 & -1 \\
 -1 & -1 & \frac{1}{2} \\
\end{array}
\right).
$$
Its three eigenvalues are 
$$
\left\{\frac{1}{4} \left(\sqrt{57}+3\right),\frac{3}{2},\frac{1}{4} \left(3-\sqrt{57}\right)\right\}
$$
and the corresponding eigenvectors are 
$$
\left(
\begin{array}{ccc}
 1 & -\frac{3-\sqrt{57}}{\sqrt{57}-9} & 1 \\
 -1 & 0 & 1 \\
 1 & -\frac{-\sqrt{57}-3}{\sqrt{57}+9} & 1 \\
\end{array}
\right).
$$

The polynomials $P_1^\pm(s)$ are given by 
$$
P_1^+(s)=3 \left(\sqrt{57}-7\right) s^2-\sqrt{57}+3, \qquad P_1^-(s)= 3 \left(\sqrt{57}+7\right) s^2-\sqrt{57}-3
$$
and their roots are real and are  algebraic numbers.\newline

The numerical study of the eigenvalues and eigenvectors of $\mu(n)$ indicates that the largest and smallest eigenvalues are simple for finite $n$, and that in the limit when $n\to \infty$ the following occurs:
\begin{fact} $(i)$~The functions realizing the maximum of the matrix $\mu$ are 
$$f^+(x)=\sin(\pi x), \ f^-(x)=\cos(\pi x), \ \forall x\in [0,1].$$
	$(ii)$~The functions realizing the minimum of the matrix $\mu$ are 
$$g^+(x)=\sin(3\pi x), \ g^-(x)=\cos(3\pi x), \ \forall x\in [0,1].$$
\end{fact}
One can then deduce the relevant properties of the Fourier transforms as follows. 
\begin{prop} $(i)$~The Fourier transforms $h^\pm$ of the functions $f^\pm(x+\frac 12)$ are given by 
$$
h^+(s)=\frac{2 \pi\,  \cos \left(\frac{s}{2}\right)}{\pi ^2-s^2}, \ \ h^-(s)=\frac{2\, i\, s \,\cos \left(\frac{s}{2}\right)}{\pi ^2-s^2}.
$$
$(ii)$~The Fourier transforms $k^\pm$ of the functions $g^\pm(x+\frac 12)$ are given by 
$$
k^+(y)=\frac{6 \pi  \cos \left(\frac{y}{2}\right)}{9 \pi ^2-y^2}, \ \ k^-(y)=-\frac{2 \,i \,y \,\cos \left(\frac{y}{2}\right)}{y^2-9 \pi ^2}.
$$
$(iii)$~The functions $h^\pm,k^\pm$ are entire functions all of whose zeros are real, and given by all odd multiples  of $\pi$ except $\pm \pi$ for $h^+$, $\pm 3\pi$ for $k^+$ and including $0$ for $h^-$ and $k^-$.	\newline
$(iv)$~The maximal and minimal eigenvalues of the matrix $\mu$ are $\frac 83$ and $-\frac 85$.
\end{prop}
\section{Explicit checks for $N=1,2$}
We  give concrete proofs of Theorem \ref{finmain} in the simplest cases $N=1,2$.
\subsection{Case $N=1$} For $N=1$ it is enough to treat the following matrix $ M(c)$ for $c\in \R$,
\begin{equation}\label{mc}
 M(c)=\left(
\begin{array}{ccc}
 0 & -1 & -1 \\
 -1 & c & -1 \\
 -1 & -1 & 0 \\
\end{array}
\right).
\end{equation}
 The three eigenvalues are 
$$
\left\{1,\frac{1}{2} \left(-\sqrt{c^2+2 c+9}+c-1\right),\frac{1}{2} \left(\sqrt{c^2+2 c+9}+c-1\right)\right\}
$$
and their dependence on $c$ is as follows:
\begin{center}
\includegraphics[scale=0.5]{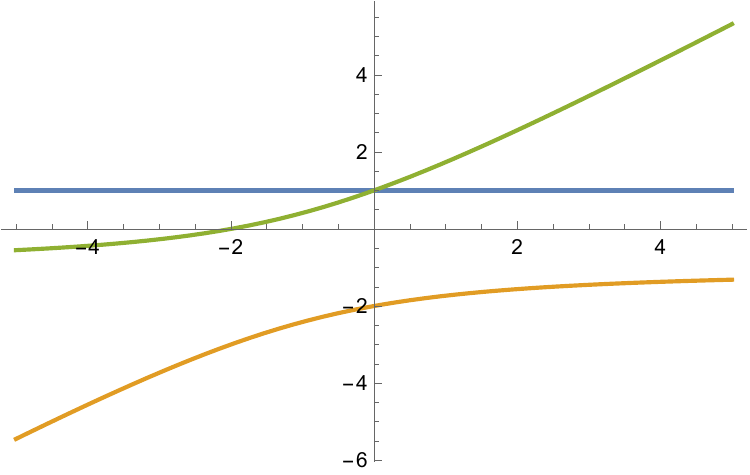}
\end{center}
The corresponding eigenvectors are 
$$
\left(
\begin{array}{ccc}
 -1 & 0 & 1 \\
 1 & X(c) & 1 \\
 1 & Y(c) & 1 \\
\end{array}
\right), \ \ X(c)=-\frac{-\sqrt{c^2+2 c+9}+c-3}{\sqrt{c^2+2 c+9}+c+3}, \ Y(c)=-\frac{-\sqrt{c^2+2 c+9}-c+3}{\sqrt{c^2+2 c+9}-c-3}
$$
and the dependence of $X(c),Y(c)$ on $c$ is as follows:
\begin{center}
\includegraphics[scale=0.5]{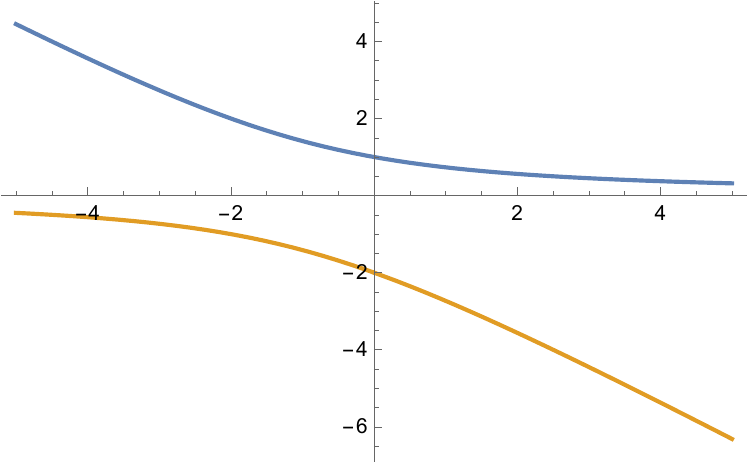}
\end{center}
The Fourier transform of the vector $(1,x,1)$ is given by the function of $s$ equal to $s^2 (-x)-2 s^2+4 \pi ^2 x$, and thus its roots are real exactly when $x(x+2)\geq 0$. This fails when $-2<x<0$, and this is realized by $Y(c)$ for $c<0$. But the corresponding eigenvalue $y(c)$ fulfills for $c<0$ 
$$
\frac{1}{2} \left(-\sqrt{c^2+2 c+9}+c-1\right)=x(c)<y(c)=\frac{1}{2} \left(\sqrt{c^2+2 c+9}+c-1\right)< 1
$$
as shown by the graphs for $c<0$, where $x(c)$ appears in blue, and $y(c)$ remains between $x(c)$ and $1$.
\begin{center}
\includegraphics[scale=0.5]{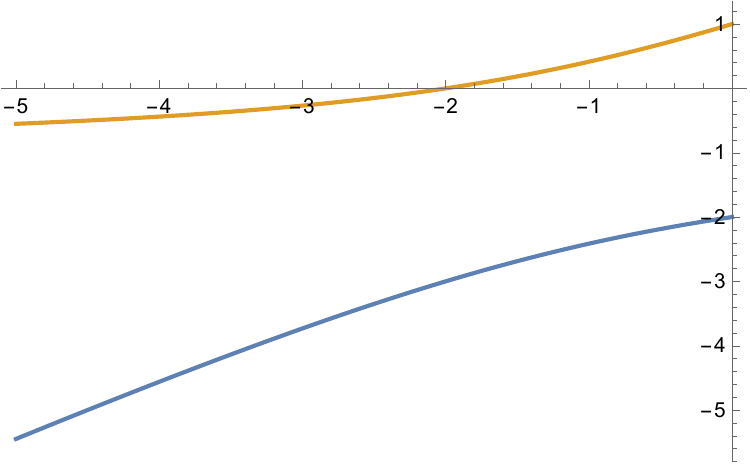}
\end{center}
We thus get 
\begin{prop}\label{simplest}For any distribution $\cD$ the $3 \times 3$ matrix $q_{n,m}$ where $n,m\in \{-1,0,1\}$ has the property that the zeros of the Fourier transforms of the eigenvectors corresponding to the extremal eigenvalues are real.	
\end{prop}
\proof Note first that it is enough to prove the result for the specific matrix $M(c)$ of \eqref{mc}. Indeed in general the $3 \times 3$ matrix $R=q_{n,m}$ where $n,m\in \{-1,0,1\}$ has all its off-diagonal entries given by Proposition \ref{matrixcomp}, and they only depend upon $\psi(1)=:a$ since $\psi(0)=0$ and $\psi(-1)=-\psi(1)=-a$. This shows 
$$
 \frac{\psi(n)-\psi(m)}{n-m}=a, \ \ \forall n,m\in \{-1,0,1\}, \mid n\neq m.
$$
If $a=0$ the eigenvectors corresponding to the extremal eigenvalues of $R$ are elements of the basis, and their Fourier transform is, by \eqref{four} a multiple of
 $$\frac{\sin( s/2)}{s/2-n\,\pi }, \ \ n\in \{-1,0,1\}$$
  whose zeros are all real. We can thus assume that $a=-1$. The diagonal values of $R$ are an even function of $n\in \{-1,0,1\}$ and thus by adding a scalar multiple of the identity matrix one can assume that they are of the form $\{0,c,0\}$ for some $c\in \R$. Let us now prove the result for the specific matrix $M(c)$ of \eqref{mc}. 
 The value of $X(c)$ is always positive and thus the Fourier transform of the vector $(1,X(c),1)$ has all its roots real for any value of $c$. The value of $Y(c)$ belongs to the forbidden interval $(-2,0)$ exactly when $c<0$, but in this case the corresponding eigenvalue $y(c)$ fails to be extremal since it is strictly between the two others $x(c)$ and $1$. \endproof 
 \subsection{Convexity proof of Proposition \ref{simplest}}
 We consider the positive cone $C_+$ in the linear space $C$ of matrices of the form 
 $$
 \mu(a,b,c)=\left(
\begin{array}{ccc}
 a & c & c \\
 c & b & c \\
 c & c & a \\
\end{array}
\right).
 $$
 We first rewrite this matrix in terms of the two orthogonal projections
 $$
P=\left(
\begin{array}{ccc}
 0 & 0 & 0 \\
 0 & 1 & 0 \\
 0 & 0 & 0 \\
\end{array}
\right), \ \  Q=\left(
\begin{array}{ccc}
 \frac{1}{3} & \frac{1}{3} & \frac{1}{3} \\
 \frac{1}{3} & \frac{1}{3} & \frac{1}{3} \\
 \frac{1}{3} & \frac{1}{3} & \frac{1}{3} \\
\end{array}
\right)
 $$
 and the identity matrix $1_3$. One gets 
 $$
 (b-a)P +3 c Q+(a-c)1_3=\mu(a,b,c)
 $$
 The ranges of the  projections $P,Q$ generate the two dimensional subspace $S$ which is the orthogonal of the vector $v=(1,0,-1)$ which belongs to the kernel of $P$ and $Q$.
 The angle of the two projections $P,Q$ is determined by its sine square, 
 $$
 (P-Q)^2=\left(
\begin{array}{ccc}
 \frac{1}{3} & 0 & \frac{1}{3} \\
 0 & \frac{2}{3} & 0 \\
 \frac{1}{3} & 0 & \frac{1}{3} \\
\end{array}
\right)=\frac 23 E, \qquad E=\left(
\begin{array}{ccc}
 \frac{1}{2} & 0 & \frac{1}{2} \\
 0 & 1 & 0 \\
 \frac{1}{2} & 0 & \frac{1}{2} \\
\end{array}
\right),
 $$
 where we denote by $E$ the orthogonal projection on $E$. By construction $E$ commutes with $\mu(a,b,c)$. We take the orthonormal basis of $E$ given by the vectors $v_1=(0,1,0)$ and $v_2=(\frac{1}{\sqrt{2}},0,\frac{1}{\sqrt{2}})$. In this basis the 
 projections $P,Q$ are given by the matrices 
 $$
 p=\left(
\begin{array}{cc}
 1 & 0 \\
 0 & 0 \\
\end{array}
\right), \ \ q=\left(
\begin{array}{cc}
 \frac{1}{3} & \frac{\sqrt{2}}{3} \\
 \frac{\sqrt{2}}{3} & \frac{2}{3} \\
\end{array}
\right).
 $$
 Any real symmetric  $2 \times 2$ matrix can be written as
 $$
 \left(
\begin{array}{cc}
 x & y \\
 y & z \\
\end{array}
\right)=\frac{1}{2} \left(2 x+\sqrt{2} y-2 z\right) p+\frac{3 y}{\sqrt{2}}q+ (z-\sqrt{2} y)1_2.
 $$
 \begin{lem} $(i)$~The map $\rho:C\to \operatorname{End}(E)$ given by $\rho(T):=ETE$ is an isomorphism of $C$ with the linear space $S(E)$ of selfadjoint real matrices in $E$.\newline
 $(ii)$~The image $\rho(C_+)$ is the intersection of the cone $S(E)_+$ of positive elements of $S(E)$, with the half space $H$, 
 $$
 H:=\{ \left(
\begin{array}{cc}
 x & y \\
 y & z \\
\end{array}
\right)\mid z\geq \sqrt{2} y\}.
 $$ 
 $(iii)$~The extreme rays of $C_+$ are transformed by the isomorphism $\rho$ into those extreme rays of $S(E)_+$ which are included in $H$.	
 \end{lem}
 
 \begin{center}
\includegraphics[scale=0.6]{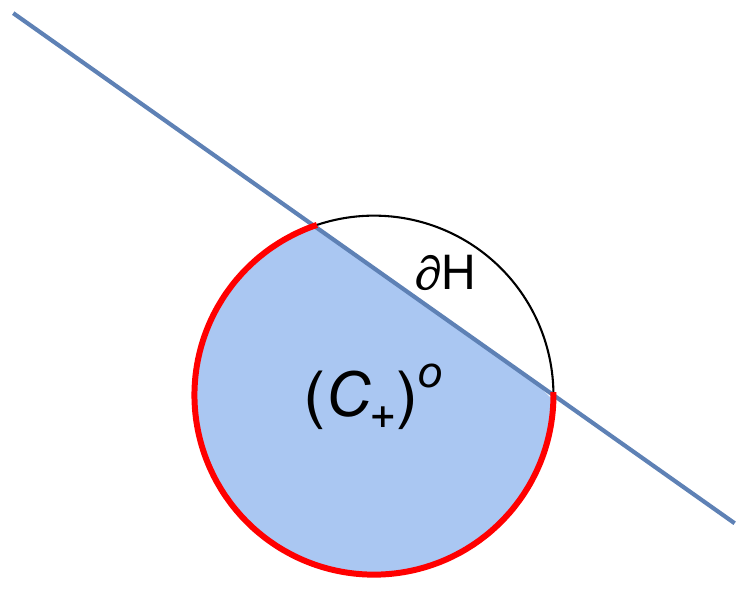}
\end{center}
\begin{prop} The elements $T=\mu(a,b,c)\in C_+$ which have a non-zero kernel are of three kinds:
\begin{enumerate}
	\item If $a>0$, then $T$ is on an extreme ray of $C_+$, its restriction to $E$ is a positive multiple of a rank one projection whose kernel coincides with the kernel of $T$.
	\item If $a=0$ and $T$ is not on an extreme ray of $C_+$, the kernel of $T$ is one dimensional, equal to $\R v$, $v=(1,0,-1)$.
	\item If $a=0$ and $T$ is  on an extreme ray of $C_+$, the kernel of $T$ is $2$-dimensional,  
\end{enumerate}	
\end{prop}
Using this proposition one obtains another proof of Proposition \ref{simplest}.
\subsection{The case $N=2$}
We are considering the linear space $C$ of real symmetric matrices of the form
$$
\mu(a,b,y,z,t)=\left(
\begin{array}{ccccc}
 t+z & b-a & \frac{b}{2} & \frac{a+b}{3} & \frac{b}{2} \\
 b-a & t+y & a & a & \frac{a+b}{3} \\
 \frac{b}{2} & a & t & a & \frac{b}{2} \\
 \frac{a+b}{3} & a & a & t+y & b-a \\
 \frac{b}{2} & \frac{a+b}{3} & \frac{b}{2} & b-a & t+z \\
\end{array}
\right).
$$
The block decomposition using the subspace $E$ for $N=1$ works in general and corresponds to the restriction to even and odd vectors, coming from the commutativity of $\mu(a,b,y,z,t)$ with the symmetry 
$$
J=\left(
\begin{array}{ccccc}
 0 & 0 & 0 & 0 & 1 \\
 0 & 0 & 0 & 1 & 0 \\
 0 & 0 & 1 & 0 & 0 \\
 0 & 1 & 0 & 0 & 0 \\
 1 & 0 & 0 & 0 & 0 \\
\end{array}
\right), \qquad J^2=\id.
$$
We take the orthonormal basis of even vectors, \ie $E:=\{\xi\mid J\xi= \xi\}$ given by 
$$
e_1=(0,0,1,0,0), \ \ e_2=(0,\frac{1}{\sqrt{2}},0,\frac{1}{\sqrt{2}},0), \ \ e_3=(\frac{1}{\sqrt{2}},0,0,0,\frac{1}{\sqrt{2}})
$$
and obtain the following matrix for the restriction of $\mu(a,b,y,z,t)$ to $E$, 
$$
\sigma(a,b,y,z,t)=\left(
\begin{array}{ccc}
 t & \sqrt{2}\, a & \frac{b}{\sqrt{2}} \\
 \sqrt{2}\, a & a+y+t & -\frac{2}{3}  (a-2 b) \\
 \frac{b}{\sqrt{2}} & -\frac{2}{3}  (a-2 b) & \frac{b}{2}+z+t \\
\end{array}
\right).
$$
For the odd vectors we take the orthonormal basis given by
$$
n_1=(\frac{1}{\sqrt{2}},0,0,0,-\frac{1}{\sqrt{2}}), \  \ \ n_2=(0,\frac{1}{\sqrt{2}},0,-\frac{1}{\sqrt{2}},0)
$$
and obtain the following matrix for the restriction of $\mu(a,b,y,z,t)$
$$
\alpha(a,b,y,z,t)=\left(
\begin{array}{cc}
 \frac{1}{2} (2 z-b)+t & -\frac{2}{3}  (2 a-b) \\
 -\frac{2}{3}  (2 a-b) & -a+t+y \\
\end{array}
\right).
$$
We now investigate positive matrices $\mu(a,b,y,z,t)$ whose kernel contains the vector $\xi=u e_1+v e_2+w e_3\in E$.
This condition specifies a $2$-dimensional subspace $K(u,v,w)$ of $C$ of the form
$$
K=\{\eta=(a,b,y,z,t)\mid y\to \frac{b w \left(3 \sqrt{2} v-8 u\right)}{6 u v}-\frac{a \left(3 \sqrt{2} u^2+3 u v-2 u w-3 \sqrt{2} v^2\right)}{3 u v},$$ $$
z\to \frac{a v \left(2 u+3 \sqrt{2} w\right)}{3 u w}-\frac{b \left(3 \sqrt{2} u^2+8 u v+3 u w-3 \sqrt{2} w^2\right)}{6 u w},t\to -\frac{\sqrt{2} a v}{u}-\frac{b w}{\sqrt{2} u}\}
$$
and one needs to decide if this subspace contains a positive element of $C$. A $2\times 2$ hermitian matrix is positive if and only if its two real eigenvalues are positive, and this is equivalent to the positivity of its trace and of its determinant. One applies this to the matrix $\alpha(\eta)$ and to the restriction of $\sigma(\eta)$ to the orthogonal of the vector $\xi$. One obtains in this manner the following $4$ "positivity conditions":
\begin{enumerate}
	\item Trace of $\alpha$ $\geq 0$ $$-\frac{3 \sqrt{2} a u+6 a v-2 a w+4 b w}{3 v}-\frac{-4 a v+3 \sqrt{2} b u+8 b v+6 b w}{6 w}\geq 0$$.
	\item Det  of $\alpha$ $\geq 0$.
	 $$\frac{-2 a^2 v \left(\sqrt{2} u+2 (v+w)\right)+a b \left(3 u^2+\sqrt{2} u (7 v+2 w)+8 v^2+6 v w-2 w^2\right)+2 b^2 w \left(\sqrt{2} u+2 (v+w)\right)}{3 v w}\geq 0$$
	\item Trace of $\sigma$ $\geq 0$ $$\frac{2 a \left(-3 u^2 w+\sqrt{2} u \left(v^2+w^2\right)-3 v^2 w\right)-b \left(3 u^2 v+4 \sqrt{2} u \left(v^2+w^2\right)+3 v w^2\right)}{3 \sqrt{2} u v w}\geq 0$$.
	\item Det of $\sigma$ $\geq 0$. $$-\frac{ \left(u^2+v^2+w^2\right) \left(2 \sqrt{2} a^2 v+a b \left(\sqrt{2} (w-4 v)-3 u\right)-2 \sqrt{2} b^2 w\right)}{3 u v w}\geq 0$$.
\end{enumerate}
On the other hand the polynomial associated to  $\xi=u e_1+v e_2+w e_3\in E$ is 
$$
P(s)=s^4 u+\sqrt{2} s^4 v+\sqrt{2} s^4 w-20 \pi ^2 s^2 u-16 \sqrt{2} \pi ^2 s^2 v-4 \sqrt{2} \pi ^2 s^2 w+64 \pi ^4 u
$$
which depends only on $s^2$ and has real roots when the degree $2$ polynomial 
$$
u x^2-20 \pi ^2 u x+64 \pi ^4 u+\sqrt{2} v x^2-16 \sqrt{2} \pi ^2 v x+\sqrt{2} w x^2-4 \sqrt{2} \pi ^2 w x
$$
has positive roots, \ie equivalently that the sum and product of the roots are positive,
which gives the two conditions\footnote{We use the reciprocal polynomial} 
\begin{equation}\label{realroopol}
\frac{5 u+\sqrt{2} (4 v+w)}{ u}>0, \ \ \frac{u+\sqrt{2} v+\sqrt{2} w}{ u}>0	.
\end{equation}

In fact we can assume that $w\neq 0$ since otherwise one is reduced to  the case $n=1$. Thus we take $w=1$, and \eqref{realroopol} is reduced to the following cases :
\begin{equation}\label{realroopol1}
\left(u<0\land v<\frac{-u-\sqrt{2}}{\sqrt{2}}\right)\lor \left(0<u<3 \sqrt{2}\land v>\frac{-5 u-\sqrt{2}}{4 \sqrt{2}}\right)\lor \left(u\geq 3 \sqrt{2}\land v>\frac{-u-\sqrt{2}}{\sqrt{2}}\right).
\end{equation}
\vspace{0.5cm}
\begin{center}
\includegraphics[scale=0.6]{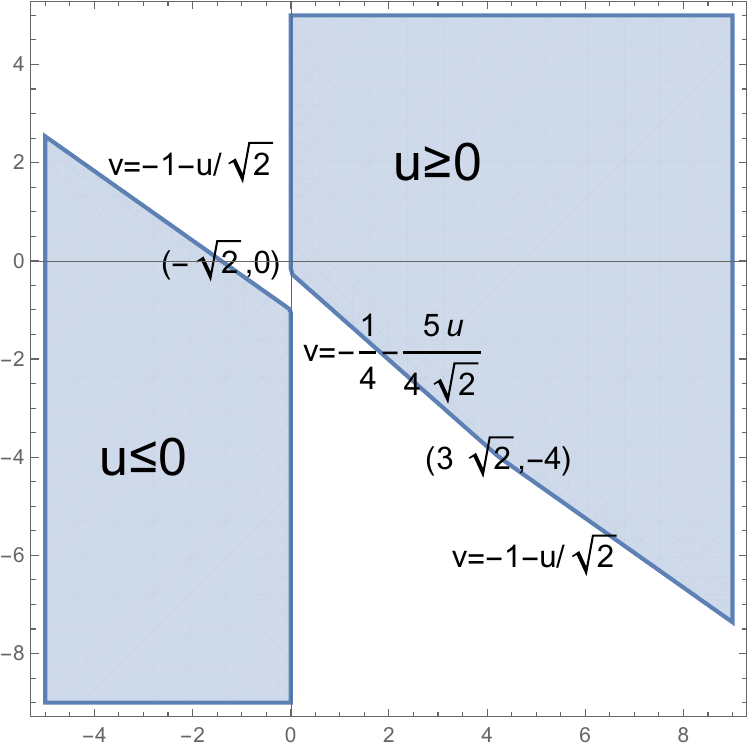}
\end{center}
This shows the region of the $(u,v)$ plane near the origin which ensures that the zeros of the associated polynomial $P$ are real.

We now reduce the above 4 positivity conditions, and compare them with \eqref{realroopol1}.
\subsubsection*{$u\leq -\sqrt{2}$}
In this case the reduction of the above 4 positivity conditions gives
$$
\left(v<0\lor 0<v<\frac{1}{2} \left(\sqrt{2} (-u)-2\right)\right)\left\|u=-\sqrt{2}\land v<0\right.
$$
which implies \eqref{realroopol1}. 
\subsubsection*{$ -\sqrt{2}<u< 0$}
In this case the reduction of the above 4 positivity conditions gives
$$
-\sqrt{2}<u<0\land v<\frac{1}{2} \left(\sqrt{2} (-u)-2\right)
$$
which implies \eqref{realroopol1}. 
\subsubsection*{$0<u<3\sqrt{2}$}
In this case the reduction of the above 4 positivity conditions gives
$$
\left(0<u<\frac{\sqrt{2}}{3}\land \left(\frac{1}{8} \left(-3 \sqrt{2} u-2\right)+\frac{\sqrt{3} \sqrt{u}}{2 \sqrt[4]{2}}<v<0\lor v>0\right)\right)\lor \left(u=\frac{\sqrt{2}}{3}\land v>0\right)$$ $$\lor \left(\frac{\sqrt{2}}{3}<u\leq 3 \sqrt{2}\land \left(\frac{1}{8} \left(-3 \sqrt{2} u-2\right)+\frac{\sqrt{3} \sqrt{u}}{2 \sqrt[4]{2}}<v<0\lor v>0\right)\right)
$$
and the solutions form the upper graph of the function $$f(u)=\frac{1}{8} \left(-3 \sqrt{2} u-2\right)+\frac{\sqrt{3} \sqrt{u}}{2 \sqrt[4]{2}}.$$
Thus to compare with \eqref{realroopol1} we need to see if the graph of $f$ in the interval $[0,3\sqrt{2}]$ is above the graph of the function  $v=-\frac{5 u}{4 \sqrt{2}}-\frac{1}{4}$. This is shown by the following
\begin{center}
\includegraphics[scale=0.6]{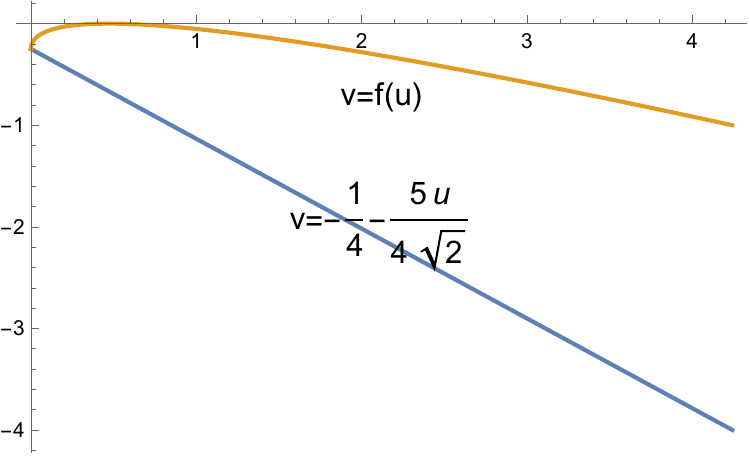}
\end{center} 
\subsubsection*{$u>3 \sqrt{2}$}
In this case the reduction of the above 4 positivity conditions gives
$$
 \left(\frac{1}{2} \left(\sqrt{2} (-u)-2\right)<v<\frac{1}{8} \left(-3 \sqrt{2} u-2\right)-\frac{\sqrt{3} \sqrt{u}}{2 \sqrt[4]{2}}\lor \frac{1}{8} \left(-3 \sqrt{2} u-2\right)+\frac{\sqrt{3} \sqrt{u}}{2 \sqrt[4]{2}}<v<0\lor v>0\right).
$$
In the first part, the condition $-1-u/\sqrt 2<v$ implies \eqref{realroopol1}. The second condition means that one is above the graph of the function $f(u)$. Thus we need to show that for $u>3 \sqrt{2}$ the graph of $f$ is above the graph of $-1-u/\sqrt 2$.
This follows from the inequality of the slopes 
$$
-3 \sqrt{2}/8>-1/\sqrt 2
$$
and the following 
\begin{center}
\includegraphics[scale=0.6]{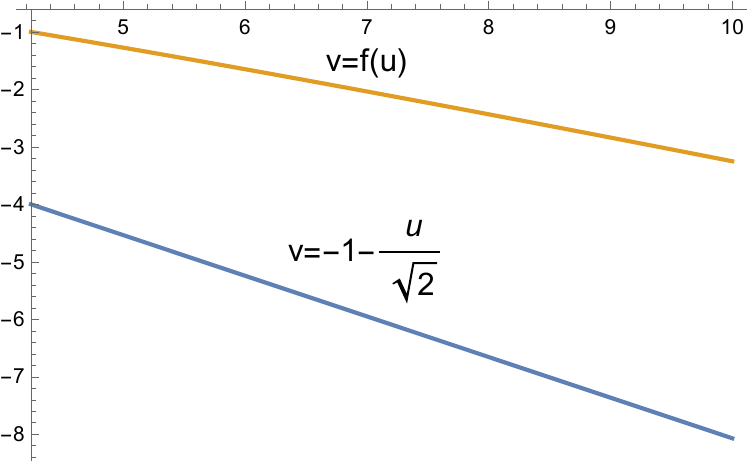}
\end{center}

We can thus conclude that if an even vector $\xi=u e_1+v e_2+w e_3\in E$ is in the kernel of a positive matrix $\mu(a,b,y,z,t)$ then the zeros of the associated polynomial are real. But we need to look at the possibility of a non trivial odd vector in the kernel of a positive matrix $\mu(a,b,y,z,t)$. \newline
We now consider the case where an odd vector $\eta=u n_1+ v n_2$ is in the kernel of a positive matrix $\mu(a,b,y,z,t)$. In this case the restriction to the odd part will be a multiple of a one dimensional projection $P_1$ and we thus need to first solve the equation 
$$
\alpha(a,b,y,z,t)=\left(
\begin{array}{cc}
 \cos ^2(\beta ) & \sin (\beta ) \cos (\beta ) \\
 \sin (\beta ) \cos (\beta ) & \sin ^2(\beta ) \\
\end{array}
\right)
$$
The solution is given by 
$$
\left\{b\to 2 a+\frac{3}{2} \sin (\beta ) \cos (\beta ),z\to \frac{1}{4} \left(-4 \sin ^2(\beta )+4 \cos ^2(\beta )+3 \sin (\beta ) \cos (\beta )\right)+y,t\to a+\sin ^2(\beta )-y\right\}
$$
We then consider the restriction of the solution matrix to the even part 
$$
\left(
\begin{array}{ccc}
 a+\sin ^2(\beta )-y & \sqrt{2} a & \frac{2 a+\frac{3}{4} \sin (2 \beta )}{\sqrt{2}} \\
 \sqrt{2} a & 2 a+\sin ^2(\beta ) & 2 a+\sin (2 \beta ) \\
 \frac{2 a+\frac{3}{4} \sin (2 \beta )}{\sqrt{2}} & 2 a+\sin (2 \beta ) & 2 a+\frac{3}{4} \sin (2 \beta )+\cos ^2(\beta ) \\
\end{array}
\right)$$
and compute its characteristic polynomial. We then apply the following
\begin{fact} Let $P(x)=x^n+\sum a_j x^{n-j}$ be a monic polynomial whose all roots are real. Then all the roots are $\geq 0$ if and only if $(-1)^ja_j\geq 0$ for all $j\in \{1, \ldots, n\}$.	
\end{fact}
One obtains in this manner the following $3$ "positivity conditions":
\begin{enumerate}
	\item  $-a_3\geq 0$ $$\frac{1}{8} \left(4 a \sin (\beta )+\cos (\beta ) \left(3 \sin ^2(\beta )-2 a\right)\right) \left(4 \sin ^3(\beta )-4 y \sin (\beta )+\cos (\beta ) \left(8 y-11 \sin ^2(\beta )\right)\right)\geq 0$$,
	\item $a_2\geq 0$
	 $$-\frac{13}{4} a \sin (2 \beta )-\left(2 a+\frac{1}{2}\right) \cos (2 \beta )-4 a y+5 a+\frac{3}{4} \sin (2 \beta )-\frac{3}{8} \sin (4 \beta )+\frac{33}{64} \cos (4 \beta )-\frac{3}{4} y \sin (2 \beta )-y-\frac{1}{64}\geq 0$$,
	\item $-a_1\geq 0$ $$5 a+2 \sin ^2(\beta )+\cos ^2(\beta )+\frac{3}{2} \sin (\beta ) \cos (\beta )-y\geq 0$$.
	\end{enumerate}
The solution of the existence of $(a,y)$ fulfilling these inequalities for a given $\beta$ is given by the following three cases. In each of them we plot the value of $-\cot (\beta )$ which gives the component $v$ of the vector $n_1+v n_2$ in the kernel of the positive matrix $\mu(a,b,y,z,t)$. We find that all values of $v$ arise except those in the interval $(-2,-\frac 12)$. The associated polynomial to the vector $u n_1+v n_2$ is $$Q(s)=2 \sqrt{2} \left(-2 \pi  s^2 u-\pi  s^2 v+8 \pi ^3 u+16 \pi ^3 v\right)$$
and its roots are real if and only if $$\frac{8 \pi ^2 (u+2 v)}{2 u+v}\geq 0$$
which for $u=1$ is realized if and only if $v\notin (-2,-\frac 12)$ as shown by the graph of the function $\frac{2 v+1}{v+2}$.

Finally  the case by case discussion of the allowed  values of $\beta$ is given below.
\begin{figure}[H]
\centering

\begin{minipage}{0.3\textwidth}
\centering
\includegraphics[width=\linewidth]{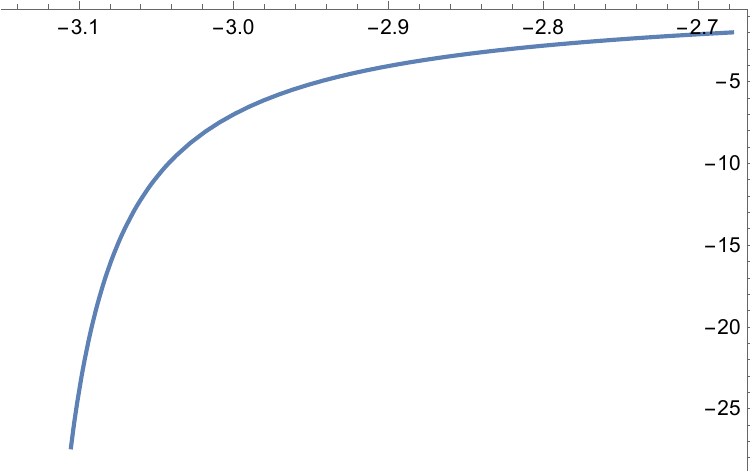}
\vspace{2mm}
\[
-\pi < \beta < -2\, \mathrm{ArcTan}\left(\sqrt{5}+2\right)
\]
\end{minipage}
\hfill
\begin{minipage}{0.3\textwidth}
\centering
\includegraphics[width=\linewidth]{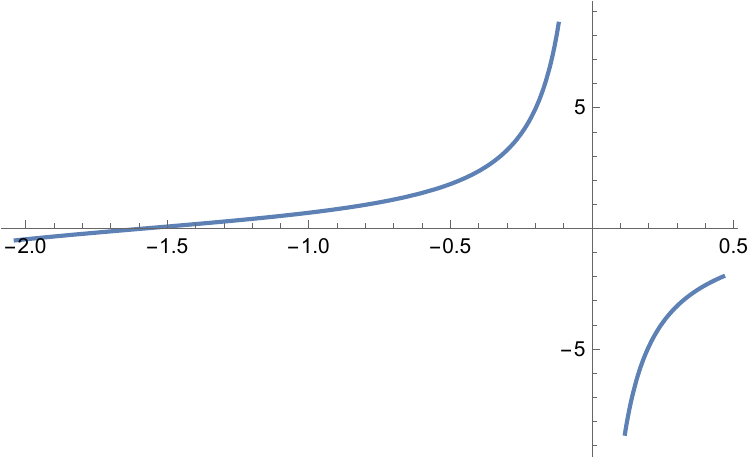}
\vspace{2mm}
$$
2\, \mathrm{ArcTan}\left(\frac{-\sqrt{5}-1}{2}\right) < \beta $$ $$ < -2\, \mathrm{ArcTan}\left(2-\sqrt{5}\right)
$$
\end{minipage}
\hfill
\begin{minipage}{0.3\textwidth}
\centering
\includegraphics[width=\linewidth]{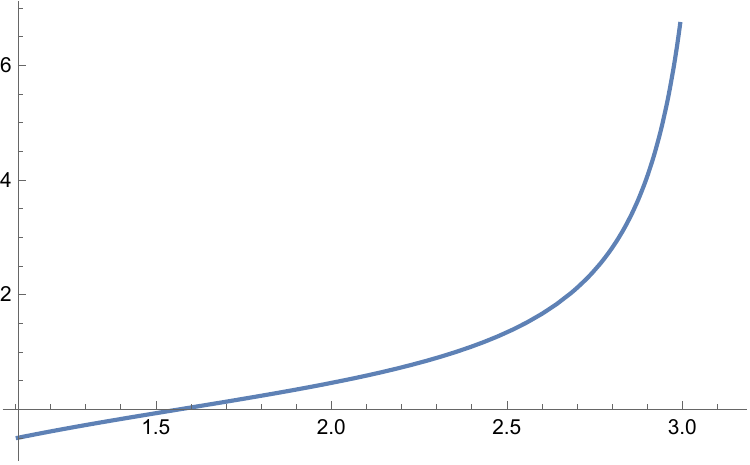}
\vspace{2mm}
$$
2\, \mathrm{ArcTan}\left(\frac{\sqrt{5}-1}{2}\right) < \beta \leq  \pi   
$$
\end{minipage}

\caption{Three conditions on $\beta$ illustrated with corresponding graphs.}
\end{figure}


\begin{thebibliography}{99}

\bibitem{araki}
  H.~Araki. Expansional in banach algebras. Annales scientifiques de l'École Normale Supérieure, Serie 4, Volume 6 (1973) 67--84.
  
  \bibitem{AM} M.~F.~Atiyah, I.~G.~MacDonald -{ Introduction to Commutative Algebra} Addison-Wesley (1969).
  
\bibitem{AS} S.~Axler, {\em  Linear Algebra Done Right}.
Fourth edition. Undergraduate Texts in Mathematics. Springer, Cham, 2024. 

\bibitem{BW} M.~Bakonyi, H.~Woerdeman, {\em Matrix Completions, Moments, and Sums of Hermitian Squares}. Princeton University Press, Princeton (2011)

  \bibitem{CF11}
C. Carathéodory und L. Fejér, {\em über den Zusammenhang der Extreme von Harmonischen Funktionen mit ihren Koeffizienten und über den Picard-Landauschen Satz}. Rend. Circ. Mat. Palermo 32 (1911)218-239. 
  


\bibitem{CC96}
A.~H. Chamseddine and A.~Connes.
\newblock Universal formula for noncommutative geometry actions: {U}nifications
  of gravity and the {S}tandard {M}odel.
\newblock {\em Phys. Rev. Lett.} 77 (1996)  4868--4871.
  
\bibitem{ckw} A.~Connes, C.~Consani, {\em Weil positivity and trace formula, the archimedean place}. Selecta Math. (N.S.) 27 (2021), no. 4, Paper No. 77.

\bibitem{CK} A.~Connes, C.~Consani,  {\em Spectral triples and $\zeta$-cycles}. Enseign. Math. 69 (2023), no. 1-2, 93-148.


 \bibitem{Mox}
A. Connes, C. Consani and H. Moscovici, Zeta zeros and prolate wave operators: semilocal adelic operators, Ann. Funct. Anal. {\bf 15} (2024), no.~4, Paper No. 87.
  
\bibitem{CW}  A. Connes and W.D. van Suijlekom, {\em Spectral truncations in noncommutative geometry and operator systems}. Comm. Math. Phys. 383 (2021), no. 3, 2021-2067. 

  
\bibitem{Don74}
W.~F. Donoghue, Jr.
\newblock {\em Monotone matrix functions and analytic continuation}.
\newblock Springer-Verlag, New York, 1974.
\newblock Die Grundlehren der mathematischen Wissenschaften, Band 207.




\bibitem{HP} E.~Hallouin, M.~Perret, {\em A unified viewpoint for upper bounds for the number of points of curves over finite fields via Euclidean geometry and semi-definite symmetric Toeplitz matrices}. Trans. Amer. Math. Soc. 372 (2019), no. 8, 5409--5451. 


\bibitem{Han06}
F.~Hansen.
\newblock Trace functions as {L}aplace transforms.
\newblock {\em J. Math. Phys.} 47 (2006)  043504, 11.

\bibitem{hermite}
C.~Hermite.
\newblock Sur la formule d'interpolation de lagrange.
\newblock {\em J. Reine Angew. Math.} 84 (1878)  70--79.


\bibitem{NSkr21}
T.~D.~H. van~Nuland and A.~Skripka.
\newblock Spectral shift for relative {S}chatten class perturbations.
\newblock {\em J. Spectr. Theory} 12 (2022)  1347--1382.

\bibitem{NS21}
T.~D.~H. van~Nuland and W.~D. van Suijlekom.
\newblock Cyclic cocycles in the spectral action.
\newblock {\em J. Noncommut. Geom.} (online 22 December 2021),
  (arXiv:2104.09899).

  \bibitem{Schmugden} K. Schmudgen, {\em Unbounded self-adjoint operators on Hilbert space}. Grad. Texts in Math. 265, Springer, Dordrecht, 2012.

\bibitem{Shannon} C.~E.~Shannon, {\em Communication in the presence of noise}, Proc. IRE 37, (1949), 10-21.


    
\bibitem{Skr13}
A.~Skripka.
\newblock Asymptotic expansions for trace functionals.
\newblock {\em J. Funct. Anal.} 266 (2014)  2845--2866.

\bibitem{Sui11}
W.~D. van Suijlekom.
\newblock Perturbations and operator trace functions.
\newblock {\em J. Funct. Anal.} 260 (2011)  2483--2496.



\end{thebibliography}
\end{document}